\documentclass{agtart_a}
\pdfoutput=1
\usepackage{graphicx}

\title{Genus generators and the positivity of the signature}

\author{A Stoimenow}
\givenname{A}
\surname{Stoimenow}
\address{Research Institute for Mathematical Sciences\\
Kyoto University\\\newline
Kyoto 606-8502\\Japan}
\email{stoimeno@kurims.kyoto-u.ac.jp}
\urladdr{http://www.kurims.kyoto-u.ac.jp/~stoimeno/}

\volumenumber{6}
\issuenumber{}
\publicationyear{2006}
\papernumber{82}
\startpage{2351}
\endpage{2393}

\doi{}
\MR{}
\Zbl{}

\keyword{signature}
\keyword{genus}
\keyword{positive knot}
\keyword{Taniyama's partial order}
\subject{primary}{msc2000}{57M25}
\subject{secondary}{msc2000}{57N70}

\received{24 June 2006}
\revised{}
\accepted{24 October 2006}
\published{13 December 2006}
\publishedonline{13 December 2006}
\proposed{}
\seconded{}
\corresponding{}
\editor{CPR}
\version{}

\arxivreference{}  %%%% not in arXiv

\AtBeginDocument{\let\bar\wbar\let\tilde\wtilde\let\hat\what}

\makeatletter
\def\cnewtheorem#1[#2]#3{\newtheorem{#1}{#3}[section]
\expandafter\let\csname c@#1\endcsname\c@lemma}
\makeatother

\newtheorem{lemma}{Lemma}[section]
\cnewtheorem{sublemma}[lemma]{Sublemma}
\cnewtheorem{theo}[lemma]{Theorem}
\cnewtheorem{theorem}[lemma]{Theorem}
\cnewtheorem{prop}[lemma]{Proposition}
\cnewtheorem{proposition}[lemma]{Proposition}

\theoremstyle{remark}
\cnewtheorem{exer}[lemma]{Exercise}
\cnewtheorem{rem}[lemma]{Remark}
\cnewtheorem{defi}[lemma]{Definition}
\cnewtheorem{conj}[lemma]{Conjecture}
\cnewtheorem{conjecture}[lemma]{Conjecture}
\cnewtheorem{corr}[lemma]{Corollary}
\makeautorefname{corr}{Corollary}
\cnewtheorem{question}[lemma]{Question}
\cnewtheorem{exam}[lemma]{Example}

\newenvironment{eqn}{\begin{equation}}{\end{equation}}

\newbox\tempboxa
\def\vcbox#1{\setbox\tempboxa=\hbox{#1}\parbox{\wd\tempboxa}{\box\tempboxa}}

\def\myfrac#1#2{\raisebox{0.2em}{\small$#1$}\!/\!\raisebox{-0.2em}{\small$#2$}}

\let\es\enspace

\AtBeginDocument{

\let\ap\alpha
\let\bt\beta

\let\gm\gamma
\let\dl\delta
\let\Dl\Delta

\let\sg\sigma
\let\Sg\Sigma

\let\sm\setminus

\let\tl\tilde

\let\x\exists
\let\vn\varnothing}

\makeop{int}

\def\len{\text{\rm len}\,}

\let\ds\displaystyle
\def\bm{\bar t_2'}
\def\bh#1{\hbox to 1ex{\hss$#1$\hss}}

\def\cK{{\cal K}}

\def\bR{{\Bbb R}}

\def\bZ{{\Bbb Z}}
\def\bN{{\Bbb N}}

\def\tD{{\tl D}}

\def\br#1{\left\lfloor#1\right\rfloor}
\def\BR#1{\left\lceil#1\right\rceil}

\def\rato#1{\hbox to #1{\rightarrowfill}}
\def\hookrato#1{\hbox to #1{$\lhook\joinrel$\rightarrowfill}}

\begin{document}

\begin{asciiabstract}
It is a conjecture that the signature of a positive link is bounded
below by an increasing function of its negated Euler characteristic.
In relation to this conjecture, we apply the generator description for
canonical genus to show that the boundedness of the genera of positive
knots with given signature can be algorithmically partially
decided. We relate this to the result that the set of knots of
canonical genus greater than or equal to n is dominated by a finite
subset of itself in the sense of Taniyama's partial order.
\end{asciiabstract}

\begin{abstract}
It is a conjecture that the signature of a positive link is bounded
below by an increasing function of its negated Euler characteristic.
In relation to this conjecture, we apply the generator description for
canonical genus to show that the boundedness of the genera of positive
knots with given signature can be algorithmically partially
decided. We relate this to the result that the set of knots of
canonical genus $\ge n$ is dominated by a finite subset of itself in
the sense of Taniyama's partial order.
\end{abstract}

\maketitle

\section{Introduction and motivation}

\subsection{The signature growth conjecture}

Most of the introductory exposition is similar, or identical to
\cite{gsig}.  A positive link (see, eg, Cromwell \cite{Cromwell2},
Ozawa \cite{Ozawa}, Yokota \cite{Yokota}) is a link which can be
represented by a positive diagram. Such links occur in several
contexts, eg, the theory of dynamical systems (Franks and Williams
\cite{WilFr,Williams}), singularity theory (A'Campo \cite{ACampo},
Boileau and Weber \cite{BoiWeb}, Rudolph \cite{Rudolph2}), and (in some
vague and yet-to-be understood way) in 4--dimensional QFTs (Kreimer
\cite{Kreimer}). They contain the class of positive braid links (see
van Buskirk \cite{Busk}, Cromwell \cite{Cromwell}), the closure links
of positive braids\footnote{\footnotesize Beware that some authors,
for example van Buskirk \cite{Busk}, confusingly call `positive links'
what we call `braid positive links' here.  Other authors call our
braid positive links `positive braids', abusing the distinction
between braids and their closures.}.  Another important subclass of
the class of positive links are the alternating ones among them, the
special alternating links (see Cromwell \cite{Cromwell2},
Nakamura \cite{Nakamura}). Such links
have been studied largely by Murasugi \cite{Murasugi3,Murasugi2}.

Knot-theoretically, one is interested how positivity can be
detected by the examination of link invariants. 
One of the most classical such invariants is the signature $\sg$.
It was studied initially by Murasugi \cite{Murasugi2}, and is defined
in terms of Seifert matrices. It thus has a natural upper bound.
For a general link $L$, let $\chi(L)$ be the maximal Euler
characteristic of a spanning Seifert surface, and $n(L)$ the number
of components. Then $\sg(L)\le n(L)-\chi(L)$, and
if $L$ has no split components bounding disconnected Seifert surfaces
(as for example for positive or alternating links), this estimate
modifies to $\sg(L)\le n'(L)-\chi(L)$, with $n'(L)$ being the number
of split components of $L$.

The positivity of the signature on positive links (or subclasses
thereof) has been a theme occurring throughout the literature
over a long period. The first result falling into this category was
established already by Murasugi in his initial study \cite{Murasugi2}.
He showed that the upper bound in terms of $\chi$ is exact for a
special alternating link \cite{Murasugi2}. This is found not to be
true for
general positive, or positive braid links by means of simple examples.

Motivated by their study in dynamical systems, in \cite{Rudolph},
Lee Rudolph showed that (non-trivial) braid positive links have
(at least) strictly positive\footnote{There is often confusion
about the choice of sign in the definition of $\sg$. Here
(following \cite{Rudolph}, rather than \cite{Traczyk} or
\cite{CochranGompf}), we use the more natural seeming convention
that positive links have positive, and not negative $\sg$.}
signature $\sg$. This result was subsequently extended to positive
links by Cochran and Gompf \cite[corollary 3.4]{CochranGompf}. A
different proof, proposed by Traczyk \cite{Traczyk}, has unfortunately
a gap and breaks down at least partly.
(It still applies for positive braid knots, the special case settled
previously by Rudolph.) Przytycki observed the result (also for almost
positive knots) to be a consequence of Taniyama's work \cite{Taniyama},
but a draft with an account on the subject was not finished. A related
proof was written down in \cite{apos}.

It is suggestive to ask how much more the signature of positive
links can grow. One should believe in an
increase of $\sg$, in the range between the maximal value in
Murasugi's result and the mere positivity property. Some evidence
suggests the following conjecture, mentioned first explicitly
in \cite{2apos}. (See \fullref{Sw} for some discussion of
this evidence.)

\begin{conjecture}\label{CJM}(Signature Growth conjecture)
\[
\liminf_{n\to\infty}\,\min\,\{\,\sg(L)\,:\,\mbox{$L$ positive link},\ 
\chi'(L)=n\,\}\,=\,\infty\,,
\]
where for a link $L$ we set $\chi'(L):=n'(L)-\chi(L)$.
\end{conjecture}

Alternatively speaking, one asks whether
\begin{eqn}\label{Sg}
\Sg_\sg\,=\,\{\,\chi'(L)\,:\,\mbox{$L$ positive link},\ 
\sg(L)=\sg\,\}
\end{eqn}
is finite for every $\sg$.

This conjecture, although suggestive, is by no means obvious,
or easily approachable. Although $\sg$ is easily calculated
for any specific link, it has turned out difficult to make
general statements about it on large link classes. This
situation is a bit opposite to $\chi$, for which much
more general formulas are available, but whose calculation
for specific links (falling outside the ``nice'' classes)
may be very complicated.

\subsection{Concordance and Bennequin's inequality}

A famous conjecture of Milnor states that for torus knots (or more
generally knots of singularities) the smooth $4$--ball genus is equal
to the genus (or unknotting number, see Boileau and Weber
\cite{BoiWeb}). This conjecture 
was settled later by gauge-theoretic work of Rudolph \cite{Rudolph3,%
Rudolph2} and Kronheimer--Mrowka (see \cite{KroMro}), which implied
the (smooth) $4$--genus version of the Bennequin inequality \cite[%
theorem 3]{Bennequin}. This inequality gives a lower bound for
the genus in terms of a braid representation of a knot or link, and
was used in his discovery of non-standard contact structures on $\bR^3$.
For a positive knot/link, the inequality estimates sharply the
genus, and hence its newer version the $4$--genus. So one obtains
explicit formulas for these invariants, and for braid positive knots/%
links from \cite{BoiWeb} also for the unknotting/unlinking number. The
(rather obvious) discussion can be found for example in \cite{Kawamura,pos}.

One of Murasugi's original results about $\sg$ is that it is a knot
concordance invariant and estimates (from below) the $4$--genus of a
knot. The signature of torus knots (and links) was found by Gordon,
Litherland and Murasugi \cite{GLM} and Hirzebruch \cite{Hirzebruch},
and fails to provide the sharp estimate desired for Milnor's
conjecture. Many more examples illustrate that the signature does not
conform to the lower bound in Bennequin's inequality. Such examples
led to the question, encountered already in Bennequin's original work,
how to modify his inequality to be applicable also to $\sg$.  A
solution was proposed in \cite{gsig}.

Recently, new signature-type concordance invariants, giving lower bounds
for the $4$--genus, were developed from Floer homology \cite{OS} and
Khovanov's homology \cite{Rasmussen} theory. Positive knots are again
intrinsically linked to these invariants, in that this time the
$4$--genus estimate is exact for such knots. (In particular, Rasmussen's
approach gives a new, combinatorial, proof of Milnor's conjecture.)

One important difference between $\sg$ and its successors is
that only former is an invariant in the topological category,
while latter apply only in the smooth category. This difference
must be emphasized in view of the growing division in methods
to study both types of concordance, where the Floer--Khovanov
homological invariants on the smooth side contrast Levine's approach
using the algebraic knot concordance group \cite{Levine} and its
recent non-abelian modifications due to Cochran--Orr--Teichner
(see \cite{Cochran}) on the topological side.

Still serious problems to understand topological concordance, and its
difference to smooth concordance, remain. Our knowledge about this
question seems to center around Freedman's result that all knots with
trivial Alexander polynomial are topologically slice. Some are known
to be not smoothly slice. The first examples were given by Andrew
Casson in the 1980s, using work of Donaldson. It was not before
these deep results that one understood topological and smooth
concordance are not the same. More such examples followed from
later work of Fintushel--Stern \cite{FintushelStern} and Rudolph
(see \cite{Rudolph,Rudolph3}), but still they remain scarce
even by now; all have trivial Alexander polynomial and rely on
Freedman's criterion. (Recently Friedl and Teichner proposed some
new candidates, with the Alexander polynomial of $6_1$, but their
smooth non-sliceness status remains unclear so far. Their simplest
good candidate for a possibly not smoothly slice knot has a diagram
with 93 crossings.) With this
state-of-the-art, one realizes to have obtained only limited 
understanding of topological concordance, and so the study of $\sg$,
a basic topological concordance invariant, gains new motivation.

\subsection{Statement of main result\label{Sw}}

In a previous paper \cite{gsig}, we settled the case of
positive braid links in \fullref{CJM}.

\begin{theo}\label{thmn}{\rm \cite{gsig}}
\[
\liminf_{n\to\infty}\,\min\,\{\,\sg(L)\,:\,\mbox{$L$ braid positive
link},\ \chi'(L)=n\,\}\,=\,\infty\,.
\]
\end{theo}

This paper is a sequel to \cite{gsig}; its motivation and problem
setting is almost equivalent; the separation was made on the one hand
for length reasons, on the other hand because the methods
applied differ somewhat. In \cite{gsig}, we used an
extension of Bennequin's inequality to $\sg$. Here we
will use another important ingredient, the generator
description for canonical genus.

To state our main result it is helpful to understand the Growth
conjecture in terms of the finiteness of the sets \eqref{Sg}.
For simplicity, consider below only knots and replace
$\chi'=1-\chi$ by the genus $g$ in \eqref{Sg}.
We start with a few remarks on known results about $\Sg_\sg$,
providing hints to the Growth conjecture. First, the positivity
result for $\sg$ means $\Sg_0=\{0\}$. Then, slightly implicitly
in \cite{Taniyama}, and later independently and explicitly in
\cite{gen2}, the result $\sg>0$ was extended
by showing that the only positive knots of $\sg=2$ are
those of genus one. Thus $\Sg_2=\{1\}$. For $\sg=4$ the situation
is not that simple. Beside genus 2, there are some positive
knots of genus 3 with $\sg=4$, and one knot of genus $4$,
$14_{45657}$ of \cite{2apos}. Calculations of \cite{2apos}
suggest that $14_{45657}$ is in fact the only positive knot 
of genus 4 with $\sg=4$, and %. This can be confirmed by
% a series of intricate diagram transformations
% and the check of a long list of special cases, similar to
% the methods given in \cite{2apos}.
for genus $g\ge 5$
indeed there seems no further such knot, that is,
apparently $\Sg_4=\{2,3,4\}$. Although this is still
difficult to check, we will resolve the problem at least
theoretically.

Our aim will be to show how one can prove, at least in theory,
that any initial number of the sets $\Sg_\sg$ is finite,
provided this is true. (Note that if $\Sg_\sg$ is infinite,
then so is $\Sg_{\sg'}$ for any $\sg'>\sg$.) Namely, we show
that there exists an algorithmically determinable collection of
knots, such that if $\Sg_\sg$ is finite, only finitely many of
the determined knots need to be checked to establish this finiteness.

\begin{theo}[Main result]\label{thalg}
For all $n>1$ there is a set $C_n$ of positive knots with two
properties:
\begin{enumerate}
\item\label{item_1}
$C_n$ is finite and algorithmically constructible.
\item\label{item_2}
For all $\sg\in 2\bN$ we have:
\[
\x \mbox{ positive knot $K$ of genus $g\ge n$ with $\sg(K)\le\sg$}\,
\iff\,
\x\ K\in C_n\,\mbox{with $\sg(K)\le\sg$}\,.
\]
\end{enumerate}
\end{theo}

To verify, using this theorem, that $\Sg_\sg$ is finite, one uses
induction on $\sg$. We know $\Sg_{\sg'}$ for $\sg'<\sg$ by induction.
Then one examines
\begin{eqn}\label{eqa}
C\,=\,C_{\textstyle\,\max \bigcup\limits_{\sg'\le \sg}\,\Sg_{\sg'}+1}\,.
\end{eqn}
If some $K\in C$ is found with $g(K)\not\in\Sg_\sg$ and $\sg(K)=\sg$,
include $g(K)$ into $\Sg_\sg$, and repeat the search for such $K$
(with the new set $C$ updated according to \eqref{eqa}).
If $\Sg_\sg$ is finite, after some iterations no $K$ will be found.

We remark that the same type of statement is true for links of
any arbitrary fixed number of components. We do not prove it,
however, in this more general form, since this generalization
does not involve significantly new arguments, and would add
only considerable technicality to the proof. We do elaborate
on the knot case, though, giving different arguments that
contribute to making $C_n$ as small as we can. We note that
in the above theorem, one can replace $\sg$ by any of the
generalized (Tristram--Levine) signatures.

In the proof of \fullref{thalg}, Hirasawa's algorithm
\cite{Hirasawa}, that lay in the center of the signature
Bennequin inequality in \cite{gsig}, finds its application
again, this time in combination with the generator
theory for diagrams of given canonical genus, initiated
in \cite{gen1}, and then developed further in \cite{STV,SV}.
Namely, we use our result of \cite{adeq} concerning the
maximal crossing number of a generating diagram of
a given genus. It improves a previous result of \cite{SV}
in this regard, and relies heavily on Hirasawa's algorithm.
The (original) Bennequin inequality also enters into the proof.

We mentioned the relation (noted by Przytycki) between the
positivity of $\sg$ and Ta\-ni\-ya\-ma's work \cite{Taniyama}. We use
this relation to bring his partial order into the context of
our arguments. It follows from our proof of \fullref{thalg}
that Taniyama's statements about the dominance of the trefoil
and $5_1$ are the first two instances of a infinite series of
such results, namely, that the set of knots of canonical
genus $\ge n$ is dominated by a finite subset of itself (\fullref{thT}). A similar outcome, for positive knots, addresses
the partial order of Cochran and Gompf \cite{CochranGompf}
(\fullref{thY}).

Further evidence for the Growth conjecture is given by the
following result on the average value of $\sg$ for given genus.

\begin{theorem}\label{tbp}{\rm \cite{adeq}}\qua
Let
\[
P_{g,n}\,:=\,\{\,K\ \mbox{positive knot, $g(K)=g$, $c(K)\le n$}\ \}\,,
\]
where $c(K)$ denotes the crossing number of $K$. Then
\[
\lim_{n\to\infty}\,
\frac{1}{\big|P_{g,n}\big|}\,
\left(\sum_{K\in P_{g,n}}\,\sg(K)\right)
\es =\es 2g\,.
\]
\end{theorem}

(Note that $P_{g,n}$ is always finite, and becomes non-empty 
for fixed $g$ when $n$ is large enough. Note also that in
general the crossing number $c(K)$ of a positive knot $K$ may not
be admitted by a positive diagram, as shown in \cite{posex_bcr}.)

This theorem means that generically the value of
$\sg$ for fixed genus is the maximal possible. From this
point of view, the philosophy behind the Growth conjecture
is that `when the generic value is the maximal possible,
the minimal value should not be too small.' \fullref{tbp} is a consequence of a (largely unrelated to
the subject of this paper) extension of the
asymptotical denseness result for special alternating knots
in \cite{SV}, which is proved in a separate paper \cite{adeq}.

\subsection{Overview of the proof}

Consider the Growth conjecture in what follows for knots.
In this subsection, before we get into considerable technicalities,
we will give a summary of the difficulty in, and strategy for
the proof of \fullref{thalg}. (A few technical terms
occurring will be explained in \fullref{Sfr}.)

For the proof we need a method to evaluate
\begin{eqn}\label{tty}
\min\,\left\{\,\sg(K)\,:\,\mbox{$K$ positive},\ g(K)\ge n\ \right\}\,.
\end{eqn}
To do so, first we apply generator theory for canonical genus
\cite{gen1,STV,SV}. This allows to calculate
\[
\min\,\left\{\,\sg(K)\,:\,\mbox{$K$ positive},\ g(K)=g\ \right\}
\]
for any given $g$ by verifying $\sg$ on finitely many knots.
These knots, the ``generators'' of \cite{gen1}, can be
algorithmically constructed. We will not get into details about this
procedure here, since we discussed it extensively elsewhere. Briefly,
there are three methods (in increasing order of efficiency):
selecting knots from the alternating knot tables (more efficient
for small crossing number), using maximal Wicks forms
(as explained in \cite{STV}, more efficient for high
crossing number), and using thickenings of trivalent graphs and
Hirasawa's algorithm \cite{gener}.

Then we need to obtain an upper bound on the genus
necessary to check for \eqref{tty}, which will be our main effort.
It is possible that always
\begin{align}\label{ttz}
\min\,\left\{\,\sg(K)\,:\,\mbox{$K$ positive},\ g(K)=g+1\ \right\}\es&\ge
\\
\min\,\{\,\sg(K)\,:\,&\mbox{$K$ positive},\ g(K)=g\ \}\,,
\nonumber
\end{align}
and then we would need to check just genus $n$. Unfortunately,
we do not know how to prove \eqref{ttz} (or whether it is always
true). If we add a `$-2$'
on the r.h.s., then the inequality follows easily from the fact
that any genus $g+1$ diagram can be turned into a genus $g$
diagram by smoothing out a (proper) pair of crossings.
However, the negative correction term continuously ruins
the estimate with increasing genus, and thus makes it
useless with regard to the Growth conjecture. We will thus
be forced to avoid smoothings and work only with crossing
changes.

Thus we need a lower bound on the
genus of a positive diagram one can obtain by crossing
switches from a given one. To find such a bound is
considerably more difficult. We will show that
the genus decreases at most by a linear factor, which
we will be able to drop to {\footnotesize $\ds\frac{221}{41}$}$
\approx 5.39$ (\fullref{thgdec}). While for the mere
existence of such a constant a part of the proofs can be
simplified, even the value we attained with the extra effort
is still too large to make our result practical. Still
one can take practical advantage of the arguments we apply,
and we will attempt to settle the problem $\Sg_4=\{2,3,4\}$
at a later stage along these lines.

In the proof of the genus decreasing bound generator theory
finds again its application. The first step of this proof 
uses an improvement of the estimate of the maximal crossing
number of a generating diagram of given genus \cite{SV}, given
in \cite{adeq}. In the next section we review the necessary tools.

Then we need an estimate of the minimal length of simplifying
bridges/tunnels. What we like is to choose a certain piece of the strand
in a positive diagram, to switch properly crossings on it so that
it becomes a bridge or tunnel, and then shrink the bridge or tunnel
by a wave move. In case the diagram has a clasp this is trivial,
so all difficulty comes from the diagrams that have no such clasps.

Since we want after the shrinking the diagram to be positive, we
must ensure the shrunk bridge/tun\-nel to have length 1, ie, only
one intersection with the rest of the diagram. In that case with
the option between bridge and tunnel we can always achieve
the new diagram to be positive. 
\vspace{-4pt}

So we must investigate how long a piece of the strand needs
to be, so that we can shrink it to a one-crossing bridge/tun\-nel.
This is very similar to what was done by myself and
M Kidwell in \cite{brlen}. There he proved, improving
my original result, that if a knot diagram contains
a bridge of length more than $\myfrac{1}{3}$ of its
crossing number, this bridge can be rerouted (wave moved)
to a smaller one. Here, however, we must work
harder, because we want to know when such a bridge can
be rerouted to one of a single crossing. What we basically
show is that to have this much stronger condition, we
can replace Kidwell's constant $\myfrac{1}{3}$ by
$\myfrac{36}{41}$. This is the content of the Curve length
lemma in \fullref{Sft}.
\vspace{-4pt}

Note that we could also seek a bridge that can be rerouted
to length 0, but we will see
that even finding (a moderately short) one shrinkable to
length 1 is difficult
enough. If a bridge/tunnel shrinks to more than one crossing,
to ensure one can keep all crossings positive,
it is necessary to take account
on the orientations of the strands intersecting the shortened
bridge/tunnel; this seems virtually unfeasible, though. A bonus of using
a bridge/tunnel of length 1 is also that we can apply our work
to Taniyama's partial order, see \fullref{thT}.
\vspace{-4pt}

\section{General preliminaries\label{Sfr}}
\vspace{-4pt}

Here we recall several basic facts and notations.
\vspace{-4pt}

\subsection{Miscellanea}
\vspace{-4pt}

%First, we fix some general (mathematical and linguistic) terminology.

By $\br{n}$ we will mean the greatest integer not greater than $n$.
By $\BR{n}$ we will mean the smallest integer not smaller than $n$.
\vspace{-4pt}

For a set $S$, the expressions $|S|$ and $\#S$ are equivalent and both
denote the cardinality of $S$. In the sequel the symbol '$\subset $'
denotes a not necessarily proper inclusion.
\vspace{-4pt}

`W.l.o.g.'\ abbreviates `without loss of generality' and
`r.h.s.'\ (resp.\ `l.h.s.')\ `right hand-side' (resp.\ `left hand-side').

\subsection{Link diagrams}

\begin{defi}
A crossing $p$ in a knot diagram $D$ is called \emph{reducible}
(or \emph{nugatory}) if $D$ can be represented in the form
$$
\cl{\includegraphics[scale=0.9]{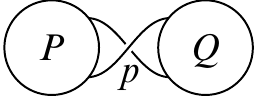}}
$$
For a nugatory crossing $p$ in $D$ there is a curve in the
plane intersecting $D$ only in $p$ (and doing so transversely).
We call this curve the \emph{nugatory curve} of $p$. A diagram $D$
is called \emph{reducible} if it has a reducible crossing, else it is
called \emph{reduced}.
\end{defi}

\begin{figure}[ht!]
\cl{\includegraphics[scale=0.9]{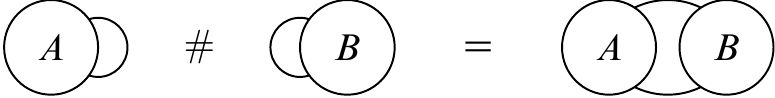}}
\caption{}\label{figtan}
\end{figure}

\begin{defi}
The diagram on the right of \fullref{figtan}
is called \emph{connected sum} $A\# B$ of the diagrams $A$ and $B$.
If a diagram $D$ can be represented as the connected sum of 
diagrams $A$ and $B$, such that both $A$ and $B$ have at least one
crossing, then $D$ is called \emph{disconnected} (or composite), else
it is called \emph{connected} (or prime).
\end{defi}

Alternatively, if $D=A\# B$ then there is a closed curve $\bt$ in the
plane intersecting $D$ in two points (and doing so transversely),
such that $A$ and $B$ are contained in the in/exterior of $\bt$.
We call $\bt$ a \emph{separating curve} for $D$.

Note in particular that prime diagrams are reduced.

\begin{defi}
If there is a closed curve $\bt$ in the plane intersecting $D$
nowhere and containing at least one component of $D$ in both its
interior and exterior, we say that $D$ is \emph{split} and $\bt$ the
\emph{splitting curve} for $D$. A \emph{split component} of a link
$L$ is a maximal set $S$ of components of $L$ with the property that
if $a,b\in S$, then in no split diagram $D$ of $L$ with $\bt$ as
splitting curve, $a$ and $b$ land on different sides in $\bR^2\sm
\bt$. A link is \emph{split} if it has a split diagram, or
equivalently, if it has more than one split component.
\end{defi}

Consider 3 links differing just at one crossing.
\begin{eqn}\label{Lpm0}
\lower 18pt\hbox{\includegraphics[scale=0.9]{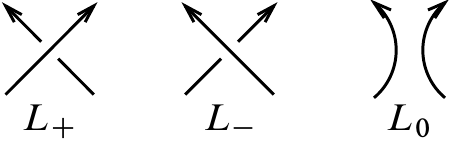}}
\end{eqn}
We call such a triple of links a \emph{skein triple}.

A \emph{positive} resp.\ \emph{negative} crossing is the fragment
of $L_+$ resp.\ $L_-$ shown in \eqref{Lpm0}. Replacing any of these
fragments by the fragment of $L_0$ in \eqref{Lpm0} is called
\emph{smoothing out} the crossing. The \emph{number of crossings}
of a diagram $D$ is written $c(D)$. The sum of signs of all
crossings in $D$ is called \emph{writhe} of $D$ and denoted
by $w(D)$.

Smoothing out all crossings in $D$ one obtains a collection of
loops in the plane called \emph{Seifert circles}. We write
$s(D)$ for the number of Seifert circles of a diagram $D$.

A knot is called
\emph{positive} if it has a diagram with all crossings positive.
Such a diagram is called also positive. A diagram is called
\emph{$n$--almost positive}, if it has exactly $n$ negative crossings.
A positive diagram $D$ obtained from a diagram $D'$ by crossing
changes is called \emph{positification} of $D'$.

A diagram is called \emph{special} if all its Seifert circles have
empty interior or exterior. Such Seifert circles in an arbitrary
diagram are called \emph{non-separating}, the others are called
\emph{separating}. Any link diagram decomposes as the \emph{Murasugi
sum} ($*$--product) of special diagrams (see \cite[Section 1]{Cromwell2}).

We call a Seifert circle $A$ \emph{opposite} to another Seifert circle
$B$ at a crossing $p$, if $p$ joins $A$ and $B$.

The \emph{(canonical) Euler characteristic} $\chi(D)$ of a link diagram
$D$ is defined as $\chi(D)=s(D)-c(D)$, where $s(D)$ is, as before, the
number of Seifert circles and $c(D)$ the number of crossings of $D$. If
$D$ is a diagram of a link with $n$ components, the \emph{(canonical)
genus} $g(D)$ of $D$ is given by 
\[
g(D)\,=\,\frac{2-n-\chi(D)}{2}\,=\,\frac{2-n+c(D)-s(D)}{2}\,.
\]
These are the genus and Euler characteristic of the 
\emph{canonical Seifert surface} of $D$, the one obtained by applying
Seifert's algorithm on $D$. The \emph{genus} $g(L)$ and \emph{Euler
characteristic} $\chi(L)$ of a link $L$ are the minimal genus and
maximal Euler characteristic of all Seifert surfaces of $L$, and the
\emph{canonical genus} $g_c(L)$ and \emph{canonical Euler characteristic}
$\chi_c(L)$ of $L$ are the minimal genus and maximal Euler
characteristic of all canonical Seifert surfaces of $L$, ie,
all Seifert surfaces obtained by applying Seifert's algorithm on
some diagram $D$ of $L$.

The importance of the canonical genus relies on the
following classical fact:

\begin{theorem}%{\rm (\cite{Crowell,Murasugi4})}
For an alternating/positive knot or link $L$ with an
alternating/positive
diagram $D$ we have $g(D)=g(L)$. (In particular, for such
knots or links canonical genus and ordinary genus coincide.)
\end{theorem}

In the alternating case this was proved by \cite{Crowell,Murasugi4}.
It can also be proved, in both cases, using \cite{Gabai}.
For positive diagrams (and in particular positive braid
representations) it follows from \cite{Cromwell}, or from
\emph{Bennequin's inequality}. The original form of this inequality
is stated as follows. 

\begin{theorem}{\rm \cite[theorem 3]{Bennequin}}\qua
If $\bt$ is a braid representation of a link $L$, then
\[
\chi(L)\le n(\bt)-|[\bt]|\,,
\]
where $n(\bt)$ is the number of strands of $\bt$ and $[\bt]$
its algebraic crossing number (exponent sum).
\end{theorem}

This inequality admits several improvements. A first, and
easy, observation is that by the braid algorithms
of Yamada \cite{Yamada} and Vogel \cite{Vogel} we obtain
a version for a general link diagram $D$ of $L$:
\begin{eqn}\label{BQ}
\chi(L)\le s(D)-|w(D)|\,.
\end{eqn}
We will use only this version of the inequality. 

Later Rudolph \cite{Rudolph3,Rudolph2} showed that the r.h.s.\
in Bennequin's inequality is actually an estimate for the
(smooth) slice Euler characteristic.
\[
\chi_s(L)\le s(D)-|w(D)|\,.
\]
This inequality was further extended by showing that one can replace
the l.h.s.\ with the invariants of Ozsv\'ath--Szab\'o and Rasmussen
on the one hand, and by slightly improving the r.h.s.\ on the other hand
(adding a strongly negative Seifert circle term; see \cite{Rudolph2,%
Kawamura}).

A \emph{clasp} is a tangle made up of two crossings. According to the
orientation of the strands we distinguish between reverse and
parallel clasp.
$$\includegraphics[scale=1.0]{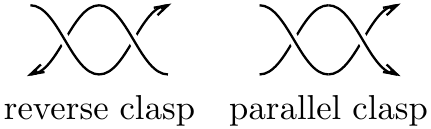}$$
By switching one of the crossings in a clasp, one can eliminate
the pair by a Reidemeister II move, and this procedure is
called \emph{resolving} the clasp.

\begin{defi}
A \emph{shadow} of a link diagram is the plane curve of
the diagram (ie, the object obtained by ignoring crossing
information). A \emph{region} of a link diagram $D$ is a
connected component of the complement of the shadow of $D$.
An \emph{edge} of $D$ is the part of the plane curve of $D$
between two crossings (clearly each edge bounds two regions).
A region is a \emph{bigon} if it has only two corners. (A bigon
in the shadow of $D$ corresponds to a clasp in $D$.)

At each crossing $p$, exactly two of the four adjacent
regions contain a part of the Seifert circles near $p$.
We call these the \emph{Seifert circle regions} of $p$. 
The other two regions are called the
\emph{non-Seifert circle regions} of $p$. 
We call two regions \emph{opposite} at a crossing $p$, if 
$p$ lies in the boundary of both regions, but they do not
share any of the four edges bounded by $p$. If two
regions share an edge, they are called \emph{neighbored}.

A diagram has a (canonical up to interchanging colors)
black-white region coloring, given by the condition
that neighbored regions have different colors. This is
called the \emph{checkerboard coloring}.
\end{defi}

A \emph{bridge/tunnel\kern0.8pt}\footnote{The term `tunnel' has usually
a different meaning in knot theory. We nevertheless use it here
in this meaning, since it is suggestive, and the other
(classical) sort of tunnels never occur in this paper.}
of a link diagram is a piece of a strand passing exclusively
through over/undercrossings. The number of such over/undercrossings
is the \emph{length} of the bridge/tunnel.

\begin{figure}[ht1]
\cl{\includegraphics[scale=0.9]{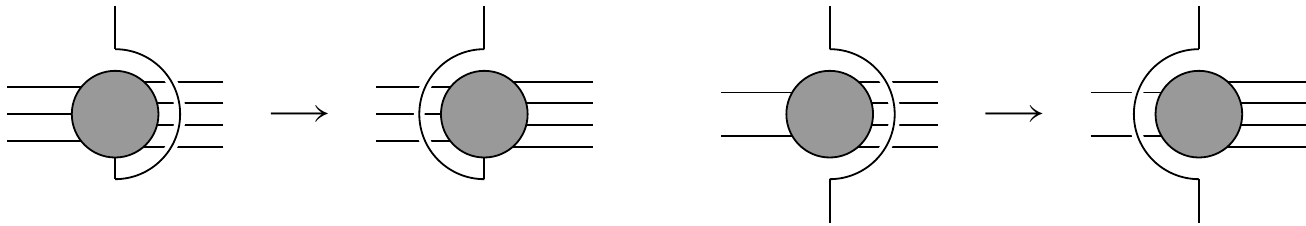}}
\caption{Wave-moves.\ The number of strands on left and right
of the shaded circle may vary. It is only important that the
parities are equal resp.\ different, and that the left-outgoing
strands are fewer that the right-outgoing ones.\label{fwv}}
\end{figure}

In certain situations, there is a move, called \emph{wave move},
that allows to shrink a bridge/tun\-nel to one of a smaller length
(see \fullref{fwv}, or also \cite{brlen} for example).

\subsection{The signature}

The \emph{signature} $\sg$ is a $\bZ$--valued invariant of knots
and links. It has several definitions. The most common one
is using Seifert surfaces and linking pairings. See, eg,
\cite{Rolfsen}. In the sequel, it will be more convenient
to follow a rather different approach, using properties of
the behaviour of $\sg$ under local (diagram) transformations.

The \emph{Alexander polynomial} $\Dl_L(t)$ can be specified by the
relation 
\[
\Dl(L_+)\,-\,\Dl(L_-)\,=\,(t^{1/2}-t^{-1/2})\,\Dl(L_0)\,,
\]
with $L_{\pm,0}$ as in \eqref{Lpm0}, and the value 1 on the
unknot. The signature $\sg(L)$ is related a value of $\Dl$,
called \emph{determinant}, $\det(L)=|\Dl_L(-1)|$. We have
that $\sg(L)$ has the opposite parity to the number of
components of a link $L$, whenever $\Dl_L(-1)\ne 0$. This
in particular always happens for $L$ being a knot ($\Dl_L(-1)$
is always odd in this case), so that $\sg$ takes only
even values on knots.

Most of the early work on the signature
was done by Murasugi \cite{Murasugi2}, who showed several properties
of this invariant.
If $L_+$, $L_-$ and $L_0$ form a skein triple, as in \eqref{Lpm0},
then\footnote{Keep in mind that our sign choice of $\sg$ follows
\cite{Rudolph} and is different from \cite{Murasugi2}.}
\begin{eqnarray}
\sg(L_+)-\sg(L_-) & \in & \{0,1,2\} \label{2a} \\[1mm]
\sg(L_\pm)-\sg(L_0) & \in & \{-1,0,1\}\,. \label{2b}
\end{eqnarray}
Further, Murasugi found the following important relation between
$\sg(K)$ and $\det(K)$ for a knot $K$.
\begin{eqn}\label{3}
\begin{array}{*3c}
\sg(K)\equiv 0\,(4) & \iff & \det(K)\equiv 1\,(4) \\[2mm]
\sg(K)\equiv 2\,(4) & \iff & \det(K)\equiv 3\,(4)
\end{array}
\end{eqn}
These conditions, together with the initial value $\sg(\bigcirc)=0$
for the unknot, and the additivity of $\sg$ under split union (denoted
by `$\sqcup$') and connected sum (denoted by `$\#$')
\[
\sg(L_1\# L_2)\,=\,\sg(L_1\sqcup L_2)\,=\,\sg(L_1)+\sg(L_2)\,,
\]
allow one to calculate $\sg$ for very many links. In particular, if
we have a sequence of knots $K_i$ 
\[
K_0\to K_1\to K_2\dots\to K_n
\]
such that $K_n$ is the unknot and $K_i$ differs from $K_{i-1}$ only by
one crossing change, then \eqref{2a} and \eqref{3} allow to calculate
inductively $\sg(K_i)$ from $\sg(K_{i+1})$, if $\det(K_{i})$ is known.

From this the following property is evident for knots, which also holds
for links: $\sg(!L)=-\sg(L)$, where $!L$ is the mirror image of $L$.

\subsection{Genus generators}

Consider the set of alternating knots $K$ of genus $g(K)=g$ and
crossing number $c(K)=n$. This set was shown to have special structure
by a theorem of \cite{gen1}, discovered independently and
simultaneously by M Brittenham \cite{Brittenham}.
In order to state this theorem, we start with some classical
definitions.

By the work of Menasco and Thistlethwaite \cite{MenThis}, alternating
knots are intimately related to a diagrammatic move called flype.

\begin{defi}
A {\em flype} is a move on a diagram shown in \fullref{fig3}.

\begin{figure}[ht!]
\cl{\includegraphics[scale=0.9]{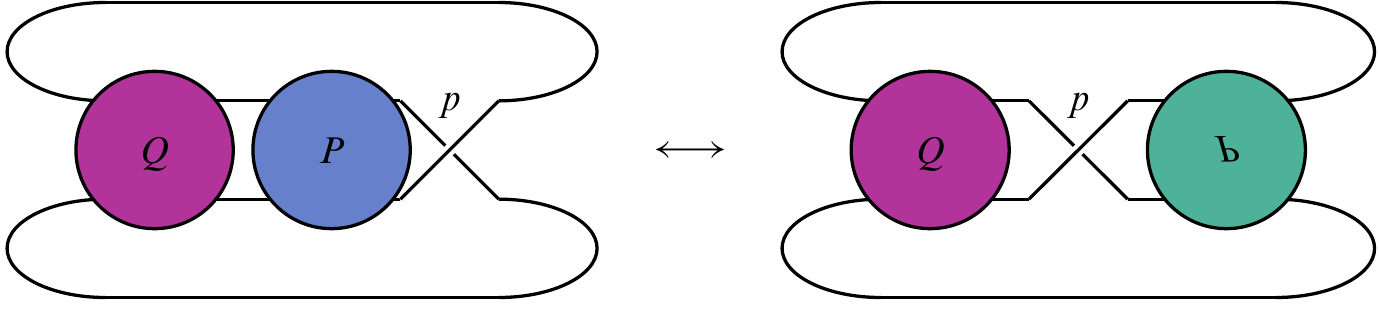}}
\caption{A flype near the crossing $p$\label{fig3}}
\end{figure}

When we want to specify the distinguished crossing $p$,
we say that it is a flype {\em near} the crossing $p$.

The tangle $P$ on \fullref{fig3} we call \emph{flypable}, and we
say that the crossing $p$ \emph{admits a flype} or that \emph{the diagram
admits a flype at (or near) $p$}. The crossing $p$ is called
\emph{flype crossing}.

\begin{defi}
A primitive Conway tangle \cite{Conway} is a tangle of the form
\[
\cl{\includegraphics[scale=1.1]{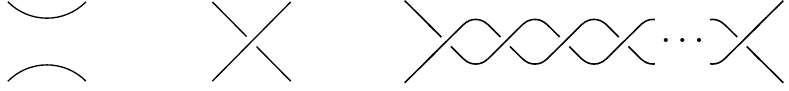}}
\]
Two crossings of a diagram are \emph{twist equivalent} if they are
contained in a primitive Conway tangle.
\end{defi}

We call the flype in \fullref{fig3} \emph{non-trivial},
if both tangles $P$ and $Q$ have 
crossings not twist equivalent to the flype crossing
(in particular they have both at least two crossings).

Since trivial flypes are of no interest we will consider
all subsequent flypes to be non-trivial, without mentioning this
explicitly each time, unless otherwise noted.
\end{defi}

\begin{theorem}{\rm \cite{gen1}}\label{th1}\qua
Reduced (that is, with no nugatory crossings)
alternating knot diagrams of given genus
decompose into finitely many equivalence classes
under flypes and (reversed) applications
of antiparallel twists at a crossing:
\begin{equation}\label{move}
\lower 10pt\hbox{\includegraphics[scale=0.9]{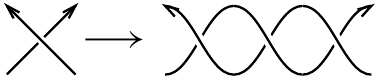}}
\end{equation}
\end{theorem}

Henceforth we call the move in \eqref{move} a \emph{$\bm$ move}.

It was observed in \cite{gen1} that in a sequence of
flypes and $\bm$ moves, all the flypes can be performed in
the beginning. It follows then from \cite{MenThis} that there
are only finitely many alternating knots with $\bm$--irreducible
diagrams of given genus $g$, and we call all such knots, and
their alternating diagrams \emph{generators} or \emph{generating}
knots/diagrams of genus $g$. The positifications of generating
diagrams of genus $g$ are called \emph{positive generating diagrams}.

There is an obvious bijective correspondence between the
crossings of the 2 diagrams in \fullref{fig3} before and after the flype,
and under this correspondence we can speak of what is
a specific crossing after the flype. In this sense, we 
make the following definition:

\begin{defi}
We call two crossings in a diagram $\sim$--{\em equivalent}, if they
can be made to form a reverse clasp after some (sequence of) flypes.
\end{defi}

It is an easy exercise to check that $\sim$ is an equivalence
relation, and that two crossings are $\sim$--equivalent iff
they share the same pair of non-Seifert circle regions.

\begin{defi}
We call an alternating diagram \emph{generating}, or a \emph{generator},
if each $\sim$ equivalence class of its crossings has $1$ or $2$
elements. The set of diagrams which can be obtained by applying flypes
and $\bm$ moves on a generating diagram $D$ we call \emph{(generating)
series} of $D$.
\end{defi}

Thus \fullref{th1} says that alternating diagrams of given
genus decompose into finitely many generating series.

\begin{defi}\label{df1}
Let $c_g$ be the maximal crossing number of a generating diagram
of genus $g$, and $d_g$ the maximal number of $\sim$--equivalence
classes of such a diagram.
\end{defi}

\section{The curve length lemma\label{Sft}}

In the following we will sometimes for convenience identify a diagram
and its shadow (its plane curve with transverse self-intersections).

Let $\gm$ be a curve, a piece of the solid line of a prime 
diagram shadow $D$. Curves are considered up to homotopy
preserving the order of edges they pass, and transversality of
their intersections. The complement $D\sm \gm$ of a curve is
a planar graph with 4--valent vertices, with two exceptional vertices
being of valence $1$. Let $v_{1,2}$ be those vertices, the start
and end of $\gm$ (we think of $\gm$ as going from $v_1$ to $v_2$).
For $D\sm\gm$ regions and edges can be defined
analogously as for $D$.

Let the \emph{length} $\len \gm$ of a curve $\gm$
be the number of its intersections with the rest of the diagram, or
the number of regions it passes (start and end region included, and
possibly reentered regions counted multiply) minus one.

\begin{defi}
Let start and end $v_{1,2}$ of $\gm$ lie in neighbored regions in
$D\sm\gm$ (that is, $\gm$ can be rerouted, or wave-moved, to a curve
of one crossing), and $\len\gm>1$. Then we call $\gm$ \emph{admissible}.
If $\gm$ has minimal length among all admissible curves we call
$\gm$ \emph{minimal admissible}.
\end{defi}

We assume for the rest of this section, unless noted otherwise, that
\[
\mbox{\emph{$\gm$ is a minimal admissible curve} and 
$D$ has \emph{no bigon regions} (clasps) \emph{and is prime}.}
\]
Under these conditions we have several lemmas.
% Let $D_{1,2}$ be the regions of $D\sm\gm$ containing the start resp.\
% end of $\gm$. We assume further that

\begin{lemma}\label{Lm1}
A minimal admissible curve $\gm$ passes any region $X$ of
$D\sm\gm$ at most once (ie, $X\cap\gm$ is empty or a single arc).
\end{lemma}

\proof Assume $\gm$ passes a region $X$ of $D\sm\gm$ twice,
and let $x_{1,2}$ be points in the interior of $X$ on $\gm$
in order of passing. Let $\tl\gm$ be the part of $\gm$ between
$x_{1,2}$. By the checkerboard coloring, $\tl\gm$ has even length,
and by assumption, this length is non-zero.
If $\tl\gm$ has length $2$, then either $D$ is composite,
or has a bigon region (clasp), which we excluded. So
$\tl\gm$ has length at least $4$. Then taking the part of
$\gm$ starting in the region after $x_1$ till $x_2$ gives a
shorter admissible curve $\gm$. \endproof

Apply the transformation \eqref{(4)} or \eqref{(4')}
to get a diagram shadow $D'$ with a
curve $\gm'$ starting and ending on edges of $D'$.
\begin{eqn}
\label{(4)}
\lower35pt\hbox{\includegraphics[scale=0.9]{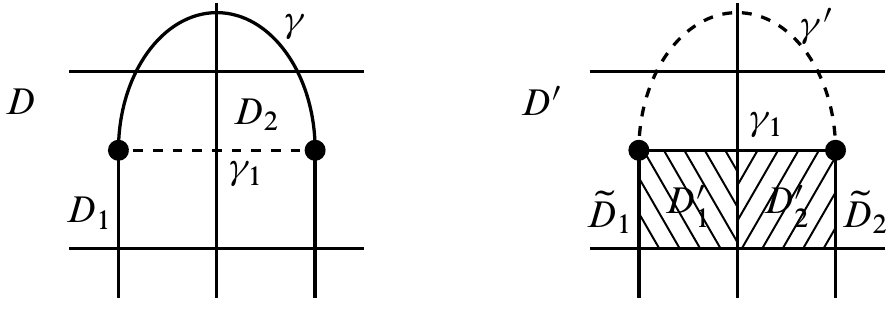}}
\end{eqn}
\begin{eqn}
\label{(4')}
\lower35pt\hbox{\includegraphics[scale=0.9]{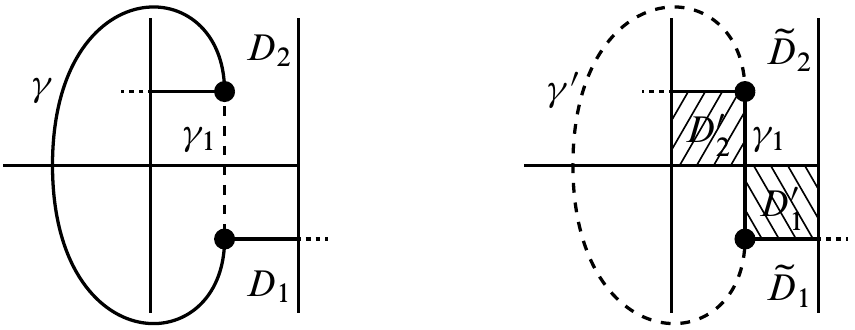}}
\end{eqn}
Let $\gm_1$ be a one crossing curve connecting the start and end
of $\gm$ in $D$. This curve $\gm_1$ is determined uniquely by the
edge is passes, but this edge may not be unique. Thus $\gm_1$
is not uniquely determined in general. It will be useful to 
remark that $\gm_1$ \emph{is} unique if $D'$ is prime, and otherwise
there are at most two possible $\gm_1$. If ambiguous, fix some
particular choice of $\gm_1$. Then let $D'=(D\sm\gm)\cup\gm_1$.
(We denote by $\gm_1$ the same arc in $D$ and $D'$, while
we distinguish between $\gm'$ and $\gm$; so $\gm'$ in $D'$ is the
``trace'' of $\gm$ in $D$.) We can assume w.l.o.g.\ that \emph{$\gm$
and $\gm_1$ do not intersect} (in interior points), otherwise a shorter
pair of curves will also do.

$D'$ has two specific regions $D'_{1,2}$ which can be described by
the property that they are bounded by a piece of $\gm_1$,
but do not contain a piece of $\gm'$ at its start or end.
Let $\tl D_{1,2}$ be the regions of $D'$ containing the start
and end of $\gm'$ and $D_i=D_i'\cup \tl D_i$. The $D_i$
can (and will) be regarded also as regions of $D\sm\gm$.

The following lemmas are suggestive from the diagrams in
\eqref{(4)} and \eqref{(4')}, but they are heavily applied in
the following arguments and still require a bit of proof.

\begin{lemma}\label{Lm1p}
The crossing on $\gm_1$ in $D'$ is not nugatory.
In other words, the four regions $D_i',\tl D_i$ for $i=1,2$
in $D'$ are all pairwise distinct.
\end{lemma}

\proof We prove indirectly.
Let this crossing, call it $x$, be nugatory, and $\dl'$
be a nugatory curve for $x$ in $D'$, and $\dl$ its preimage in
$D$. (We write $x$ also for the preimage of $x$ in $D$.) Since
$\dl'$ intersects $D'$ only in $x$, the start and end of
$\gm_1$ must lie in different regions of $\bR^2\sm\dl$.
Therefore, $|\dl'\cap \gm'|=|\dl\cap\gm|$ is odd. Now, as $\dl'$
is a nugatory curve, $\dl$ passes only one region of $D\sm\gm$.
Using \fullref{Lm1}, and a proper homotopy of $\dl'$,
we see that we can reduce the intersections between
$\dl$ and $\gm$ to one. We assumed that $\gm\cap\gm_1=\vn$,
so that $x\not\in\gm'$, and so we can choose the homotopy so that
$\dl'$ is still a nugatory curve for $x$. Now, since $\dl\cap
(D\sm\gm)=\{x\}$ and $|\dl\cap\gm|=1$, we have $|\dl\cap D|=2$.
By primality of $D$, (after pushing $\dl$ off $x$ slightly)
one of the regions of $\bR^2\sm\dl$ contains just a trivial
arc of $D$. But this arc must then contain one of the endpoints
of $\gm$ and $\gm_1$. Also, the strand of $D\sm\gm$ that intersects
$\dl$ in $x$ must continue into this arc.
So $x$ must be an intersection of $\gm_1$
with $\gm$, which we excluded. This is a contradiction. \endproof

\begin{lemma}\label{Lm2}
$\gm'$ does not enter into $D_1'$ or $D_2'$, ie,
$\gm'\cap(\int\,D_1'\,\cup\, \int\,D_2')=\vn$. 
(Here `$\int$' stands for the topological interior.)
\end{lemma}

\proof By construction the initial and terminal parts of $\gm'$
do not enter into $D_1'$. Also $\gm'\cap\gm_1=\vn$ by assumption.
Now let $e_1$ be the edge of the shadow of $D'\sm\gm_1$ that
contains the start point $v_1$ of $\gm'$ (and $\gm_1$). If $\gm'$
intersects $e_1$ in $D'$ before leaving $D_1=\tl D_1\cup D_1'$, then
$D$ is not prime, which we excluded. Since $\gm'\cap\gm_1=\vn$, 
this means that $\gm'$ must enter $D_1'$ through an edge which
lies in the boundary of $D_1$ in $D$. Then $\gm$ passes
$D_1$ in $D\sm\gm$ twice, in contradiction to \fullref{Lm1}.
The case of $D_2'$ is similar. \endproof

\begin{lemma}\label{Lm3}
The only neighbored regions to $D_{1,2}'$ passed by $\gm'$ are
$\tl D_{1,2}$.
\end{lemma}

\proof Let, contrarily, w.l.o.g.\
$X$ be a neighbored region to $D_1'$ in $D'$, different
from $\tl D_{1,2}$, that contains a part of $\gm'$. Then $X$ does not
contain an initial or terminal arc of $\gm'$. So taking the part
$\hat\gm'$ of $\gm'$ from $v_1$ to (an interior point of) $X$
would give a shorter admissible curve $\hat\gm$ in $D$, unless
$\len_D\hat\gm=1$.

Assume $\len_D\hat\gm=1$. Since $\hat\gm\cap\gm_1\subset \gm\cap\gm_1=
\vn$, then also $\len_{D'}\hat\gm'=1$. Now the curve $\hat\gm'$ starts
at $v_1$ in $D'$, and then enters either $D_1'$ or $\tD_1$. Since
$\tD_1$ has the same checkerboard color as $X$ in $D'$, the
curve $\hat\gm'$ needs at least two edges to intersect to
reach $X$ from $\tD_1$. Therefore, if $\len_{D'}\hat\gm'=1$, then
$\hat\gm'$ enters $D_1'$. Then the same is true, however, for
$\gm'\supset \hat\gm'$, but this contradicts \fullref{Lm2}.

This contradiction
shows that $\len_D\hat\gm=\len_{D'}\hat\gm'>1$. So $\hat\gm$ is
admissible in $D$, and again we have a contradiction to the minimality
of the admissible curve $\gm$. \endproof

\begin{lemma}\label{Lm3.4'}
There exists a unique edge of $D'$, the one containing $v_j$,
that bounds (simultaneously) $\tD_j$ and $D_j'$ (for $j=1,2$).
If $\tD_j$ and $D_{3-j}'$ are neighbored in $D'$, then, again,
their commonly bounding edge is unique.
\end{lemma}

\proof Consider w.l.o.g.\ only $j=1$. Let $g_i$ be the parts of $\gm_1$
in $D'$ that bound $D_i'$ and $\tD_i$ for $i=1,2$.  Let $e_i$ be the
edge of $D'\sm\gm_1$ that contains $v_i$ (or bounds $D_i$ on either
side). Clearly $x_1=g_1\cup e_1$ is an edge in $D'$ that bounds
$D_1'$ and $\tD_1$. Assume there is another such edge $x_1'$.
Then ($D'$ is composite and) there is a separating curve $\dl$ in
$D'$ that intersects $D'$ only in $x_1$ and $x_1'$. By homotopy
we can assume w.l.o.g.\ that $\dl$ intersects $x_1=g_1\cup e_1$
in some point on $e_1$, so $\dl\cap g_1=\vn$. Now, since $\dl$
passes only the regions $D_1'$ and $\tD_1$, by \fullref{Lm1p},
$\dl$ cannot intersect $g_2$. So $\dl\cap(g_1\cup g_2)=\dl\cap\gm_1=
\vn$. Then $v_{1,2}$ lie in the same region of $\bR^2\sm\dl$.

Now consider $\dl\cap\gm'$. Since $\dl$ passes only $D_1'$ and
$\tD_1$, but $\gm'$ does not enter $D_1'$,
we have $\dl\cap \gm'\subset \tD_1$. Now both $\gm'$ and $\dl$
have only one arc in $\tD_1$ (for $\gm'$ use \fullref{Lm1}),
so by homotopy we may assume that $|\dl\cap\gm'|\le 1$. But since
$\dl\cap\gm_1=\vn$ and $|\dl\cap(\gm'\cup \gm_1)|$ is even (by
Jordan curve), we see $\dl\cap(\gm'\cup \gm_1)=\vn$. Then, however,
$\dl$ remains a separating curve in $D$, and $D$ is composite,
a contradiction. 

This shows the lemma for an edge bounding $D'_1$ and $\tD_1$.
The argument for (an edge between) $D_2'$ and $\tD_1$ is similar. \endproof

The following is useful to record here, though it will be needed
only in a later stage.

\begin{lemma}\label{Lm4}
$\gm'$ does not pass a region twice in $D'$, and $\len\gm'=\len\gm$.
\end{lemma}

\proof The only difference between the regions of $D\sm\gm$
and $D'$ is that $D_i$ are subdivided into $D_i'$ and $\tl D_i$.
However, by \fullref{Lm2}, $\gm'$ does not pass $D_i'$,
so the claim follows from \fullref{Lm1}. \endproof

Let $r(D')$ be the \emph{number of regions} of $D'$. The main content of this
section, whose proof will occupy its rest, is the following estimate.

\begin{lemma}[Curve length lemma]
Assume again that $D$ is prime and that it has no clasp. Then
$\gm'$ passes 
\begin{eqn}\label{mxz}
\len \gm'+1\,\le\,\max\left(\,\frac{36r(D')-96}{41},\,4\,\right)\,
\end{eqn}
regions of $D'$, including first and last.
\end{lemma}

\proof 
We will again for convenience abuse the distinction between the
shadow $D$ and the positive diagram realizing it.
Take a prime diagram $D$, and choose $\gm$ with start and end
in neighbored regions of $D\sm\gm$, which has length $>1$.
Such a curve always exists. For example, take a loop
and consider the curve going from the loop crossing until before
the last crossing the loop passes.
\[
\cl{\includegraphics[scale=0.9]{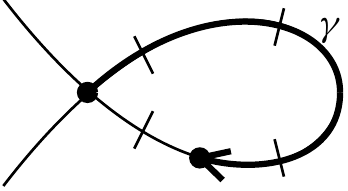}}
\]
Clearly start and end region of such a curve are
neighbored. If the curve has length $1$, then $D$ is composite
or has a clasp, in contradiction to our assumption.

We fix now among all admissible $\gm$ (not necessarily such that
come from a loop) one of minimal length ($>1$),
and the (chosen, if ambiguous) curve $\gm_1$ of one
crossing connecting start and
end of $\gm$.

When $D$ is prime, often $D'$ will also be prime.
We want to show that $D'$ becomes composite only in
a very restricted situation.

Assume that $D'$ is not prime. Let $\bt'$ be a
separating curve of $D'$. (Recall that $\bt'$ is
characterized by intersecting $D'$ in two points, and
its interior and exterior being not simple arcs.)
Thus $\bt'$ passes exactly two regions of $D'$.
If a region $X$ in $D'$ is passed by some separating curve $\bt'$,
then we call $X$ a \emph{separating region}.
\begin{eqn}\label{(8)}
\lower 35pt\hbox{\includegraphics[scale=0.9]{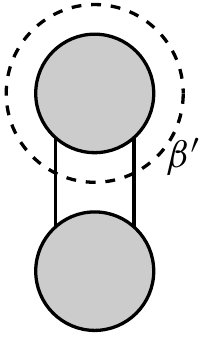}}
\end{eqn}
\begin{sublemma}\label{slp}
$D'$ has at most one pair of separating regions, and if
so, they are the start and end region of $\gm'$.
\end{sublemma}

\proof Fix a separating curve $\bt'$ in $D'$.
It has a preimage in $D$ we call $\bt$. This curve $\bt$ cannot
be a separating curve of $D$ by primeness assumption.
So the separating property must have been spoiled when recovering
$D$ from $D'$. If $\bt$ does not intersect $\gm\cup\gm_1$ in $D$,
then the move from $D$ to $D'$ must eliminate all
crossings inside or outside of $\bt$, which clearly does not happen.

% Thus $\bt$ must have in $D$ more than two crossings.
Thus $\bt$ must intersect $\gm\cup \gm_1$ in $D$.
It intersects $\gm\cup \gm_1$ in some non-zero
even number of points, %of which at least one lies in $\gm$,
and $D\cup \gm_1$ in exactly two other points not on $\gm$.

Assume now $\bt$ intersects $\gm_1$ in $D$. Then $\bt'$ does
so also in $D'$, that is, (at least) one of
the two points on $\bt'\cap D'$ lies on $\gm_1$.
\[
\cl{\includegraphics[scale=0.9]{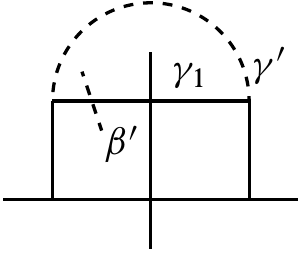}}
\]
If $\bt'$ intersects $\gm_1$ twice, then by \fullref{Lm1p}
it must do so on the same side of the crossing of $\gm_1$ in $D'$.
(Keep in mind that $\bt'$ passes only two regions in $D'$.)
Then, however, since $\bt'$ intersects $D'$ in no further points,
it cannot be a separating curve. Thus $\bt'$ intersects $\gm_1$ only
once, and $\gm'$ at least once (so that also $\bt\cap\gm\ne\vn$).

Then $\bt'$ passes through one of the regions $D'_{1,2}$,
which do not contain a part of $\gm'$. By \fullref{Lm3},
the only neighbored regions to
$D'_{1,2}$ containing a part of $\gm'$ in $D'$ are $\tl D_{1,2}$.
But if $\bt'$ passes only through $D'_i$ and $\tl D_i$ in $D'$,
then $\bt$ remains in $D$ entirely within $D_i$. Consequently,
one of the interior or exterior of $\bt'$ in $D'$ contains
only a trivial arc, a contradiction.
\[
\cl{\includegraphics[scale=0.9]{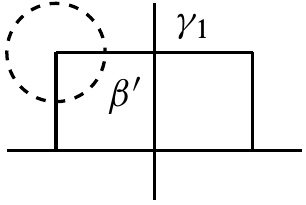}}
\]
If $\bt'$ passes through $D_i'$ and $\tl D_{3-i}$, then we could
have something like:
\[
\cl{\includegraphics[scale=0.9]{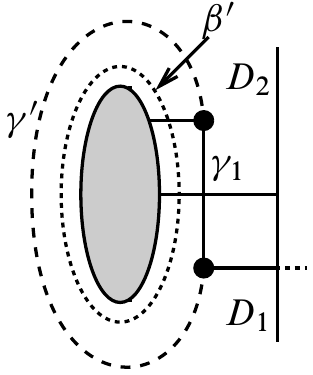}}
\]
However, we assumed that $\bt$ (and so $\bt'$) intersects $\gm_1$.
Then $\bt'$ would also pass through $\tl D_i$ or $D_{3-i}'$.
By \fullref{Lm1p}, this will be a third region
$\bt'$ passes, again a contradiction.

Therefore, $\bt$ does not intersect $\gm_1$ in $D$.
% so $\bt$
% must intersect $\gm$ (and $\bt'$ intersects $\gm'$).
Thus $\bt$ intersects non-trivially $\gm$ in $D$ 
(and $\bt'$ intersects $\gm'$ in $D'$) in an even number
of points, and two other points of $D$ not in $\gm$.
Additionally, we see that \emph{$\bt'$ does not pass $D'_{1,2}$}.
Otherwise, it would have to pass through a neighbored region
of $D'_{1,2}$ containing a part of $\gm'$. The only such regions
are $\tl D_{1,2}$, and then we have a contradiction using
\fullref{Lm3.4'} and the preceding argument.

Since $\bt'$ passes through exactly two regions in
$D'=(D\sm\gm)\cup\gm_1$, its preimage $\bt$ can pass
through at most two regions in $D\sm\gm$. If it passes
through only one region in $D\sm\gm$, then it must be
some of the $D_i$, which are subdivided in $D'$.
However, then $\bt'$ passes in $D'$ through $D'_i$,
which we argued out. Thus $\bt$ passes through exactly
two different regions $X$ and $Y$ in $D\sm\gm
$, which are therefore neighbored in $D\sm\gm$.

We claim that $\bt$ cannot intersect the two points in $D\sm\gm$
consecutively, without intersecting $\gm$ in between.
Otherwise, all intersections of $\bt$
with $\gm$ would be consecutive, say in region $X$ of $D\sm\gm$.
\[
\cl{\includegraphics[scale=0.9]{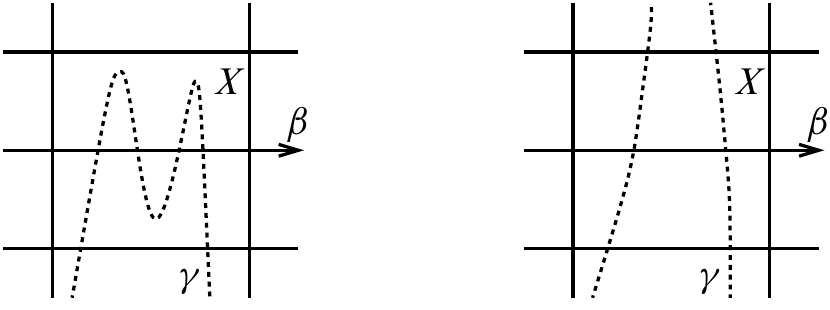}}
\]
Since by \fullref{Lm1}, $\gm$ does not reenter the
same region $X$ (ie, the second of the above two pictures does
not occur), we see (as in the first of the above pictures)
that we can homotope $\bt$ within $X$ off $\gm$, that
is, we can find a curve $\tl\bt$ in $D$ not intersecting
$\gm$, such that in $D'$ we still have the form \eqref{(8)}.
This is a contradiction to the preceding arguments that $\gm$
intersects $\bt$.

Since $\bt$ intersects $\gm$ between the two intersections with $D\sm
\gm$, this means also that both $X$ and $Y$ are passed by $\gm$.
Assume some of $X$ and $Y$ is not the start or end region of $\gm$
in $D\sm\gm$. Then by the minimality of $\gm$ and \fullref{Lm1},
$\gm$ must pass directly from $X$ to $Y$. By primeness of $D$,
there is a unique edge of $D\sm\gm$ to pass to move (directly) from
$X$ to $Y$. (See the remark below \eqref{(4')}.) So
we see that $\bt$ and $\gm$ pass from $X$
to $Y$ through the same edge of $D\sm\gm$. Since this edge is
clearly not intersected by $\gm_1$, then we can homotope $\bt'$ in $D'$
so that $\bt$ and $\gm$ do not intersect in $X$. Thus $\bt$
intersects $\gm$ only in $Y$, and then the intersections of $\bt$ with
$D\sm\gm$ are consecutive, which we argued out. 

So the two regions $X$ and $Y$ are the start and end
regions of $\gm$ in $D\sm\gm$, which we called $D_{1,2}$. When
turning $D$ into $D'$, each of $D_{1,2}$ are subdivided into
two parts by a piece of $\gm_1$. These parts are
$D'_{1,2}$ and $\tl D_{1,2}$. We argued that $\bt'$ does not pass
$D'_{1,2}$. Thus it intersects the other two parts $\tl D_{1,2}$.
By a checkerboard coloring argument, we see that we must have the
case \eqref{(4)} (rather than \eqref{(4')}), and $D$ looks like:
\[
\cl{\includegraphics[scale=0.9]{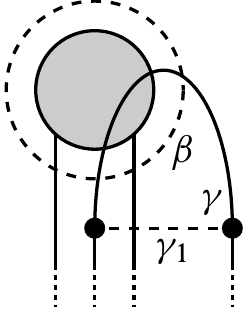}}
\]
So the start and end region $\tl D_{1,2}$ of $\gm'$ are the only
regions passed by $\bt'$. The proof of \fullref{slp} is
therefore complete. \endproof

We wish to bound from below now the number of regions of $D'$
in terms of $\len\gm'$. For this sake, we will count the 
regions of $D'$ ``close'' to $\gm'$.

\begin{figure}[ht!]
\cl{\includegraphics[scale=0.9]{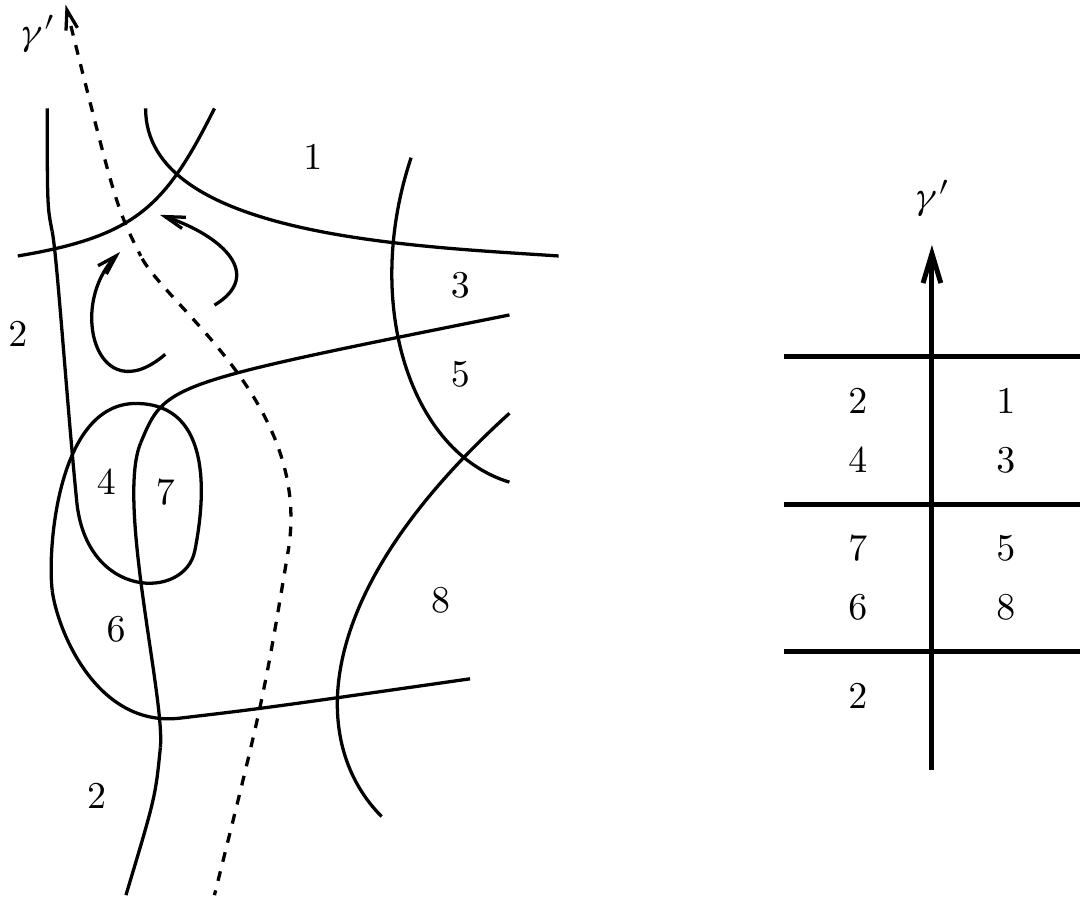}}
\caption{}\label{fig4}
\end{figure}

Number the regions of $D'$ and make a table (\fullref{fig4}). This
table contains two columns, left and right, and a row for each
intermediate region $R$ (ie, $R\ne \tD_{1,2}$) passed by $\gm'$. Note
that the number of regions $R$ is non-zero by assumption of
admissibility on $\gm$.  We call a part of the table given by the row
and choice of left/right side a \emph{slot}.

For each region $R$ passed by $\gm'$, except first and last,
order into the left and right side of the table the neighbored
regions $R'$ to $R$ from left/right of $\gm'$
in negative/positive rotation sense, as shown in the below figure.
These are the regions neighbored to $R$ by edges not intersected
by $\gm'$. The edges separating $R$ and $R'$ are unique for any pair
of neighbored regions $R$, $R'$ because we proved in \fullref{slp} that any separating curve of $D'$ does not pass through
a region $R$ passed intermediately by $\gm'$. 

In the row of the table corresponding to $R$,
put the regions on the left/right of $\gm'$ thus ordered from
bottom to top.

Then each neighbored region $R'$ entered into the table is not
passed by $\gm'$, otherwise there is a contradiction to the
minimality of its length. For the same reason, using our
preparatory lemmas, we showed that $D'_{1,2}$ are not passed
by $\gm'$, and also any neighbored region of theirs is not
passed by $\gm'$ as an intermediate region (ie, not first or
last). Thus $D'_{1,2}$ do not occur in the table.

We would like to count the regions $R'$ in the table,
but must avoid duplications. We will thus make some
effort to bound the number of such duplications. 

For this we remark that
there are several rules the entries in the table satisfy:
\begin{itemize}
\item All entries on the left and right on the same line
are distinct (because the restriction to the separating region pairs
in $D'$).
\item All entries on the left are distinct from all entries on
the right (in whatever row). This follows from the Jordan
curve theorem applied on $\gm\cup\gm_1$.
\item Each row contains at least one entry (on the left or
right), ie, there are no empty rows:
\[
\cl{\includegraphics[scale=0.9]{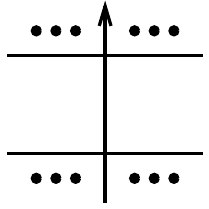}}
\]
Otherwise we have a move
\[
\cl{\includegraphics[scale=0.9]{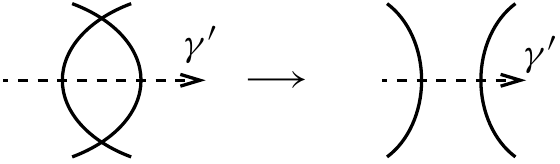}}
\]
\newpage
which shrinks an admissible curve $\gm$ of length 3 in $D$, and
since by \fullref{Lm4}, $\len\gm'=\len\gm=3$, we are done,
using the second alternative in the maximum in \eqref{mxz}. 

\item There are no two consecutive rows which are empty
on the same side, ie:
\[
\cl{\includegraphics[scale=0.9]{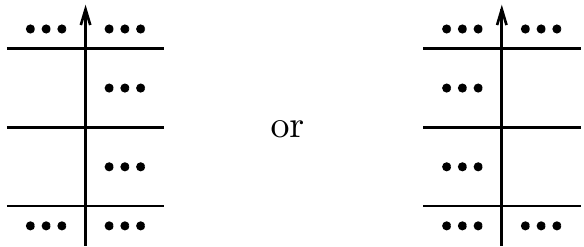}}
\]
Otherwise we can find a shorter curve $\gm$ of length 3,
which would shrink to $\gm_1$ like
\[
\cl{\includegraphics[scale=0.9]{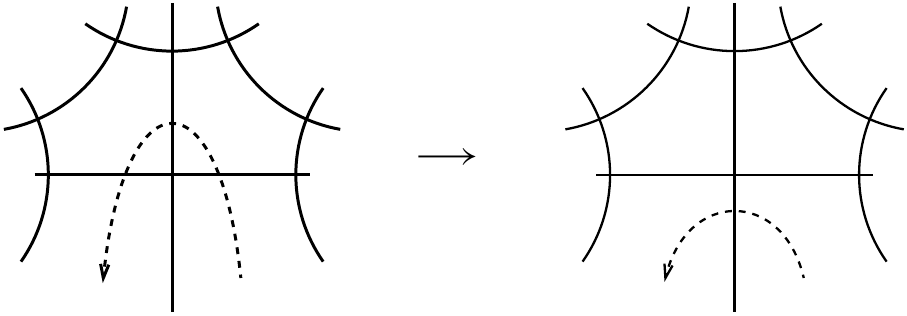}}
\]
and again are done
using the second alternative in the maximum in \eqref{mxz}. Finally,
\item there are no 3 consecutive rows with only one entry
(left or right):
\[
\cl{\includegraphics[scale=0.9]{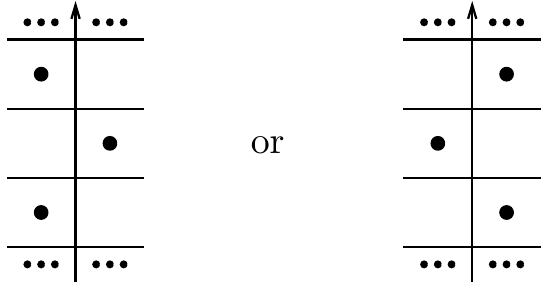}}
\]
(The previous point rules out the other patterns.)
Otherwise we have a picture like this:
\[
\cl{\includegraphics[scale=0.9]{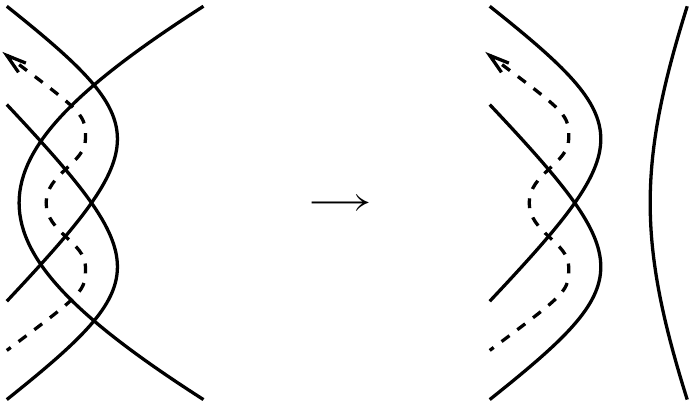}}
\]
\newpage
The right strand then contains a curve $\gm'$ (or $\gm$ in $D$)
of length $3$, and so we are done as in the previous point.
\end{itemize}

Replace now in table numbers by beads, and introduce an
equivalence relation between beads coming from the same number
($=$region). If now two beads are equivalent, then they have the
same parity of row (because of the checkerboard coloring),
and the same side (because of the Jordan
curve theorem applied on $\gm\cup\gm_1$).

Join by an arc two consecutive equivalent beads on left or
right (ie, all beads between them are not equivalent to
them).
\[
\cl{\includegraphics[scale=0.9]{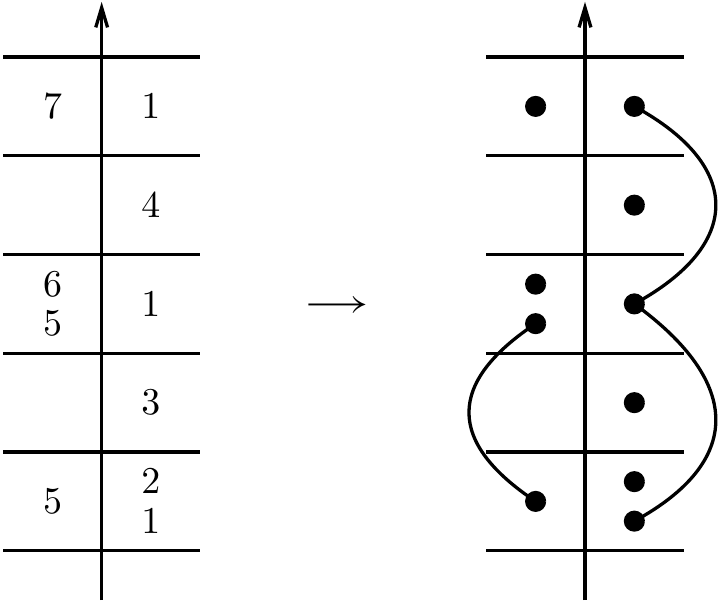}}
\]
We will consider henceforth (numerical) entries and
beads (thus connected by arcs) as equivalent, since both
representations carry the same information.

Then there are no overcrossing arcs (again by the
Jordan curve theorem, since every arc can be thought of as
lying within the same region):
\[
\cl{\includegraphics[scale=0.9]{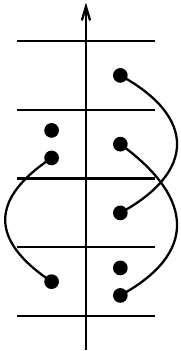}}
\]
\newpage
Finally, no two beads in the same row and side are
equivalent (because $R$ is not in a pair of separating
regions in $D'$):
\[
\cl{\includegraphics[scale=0.9]{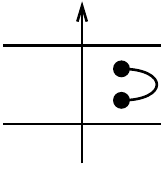}}
\]
On one of the sides the equivalence arcs may look (after
rotating by $90^\circ$) like
\[
\cl{\includegraphics[scale=0.9]{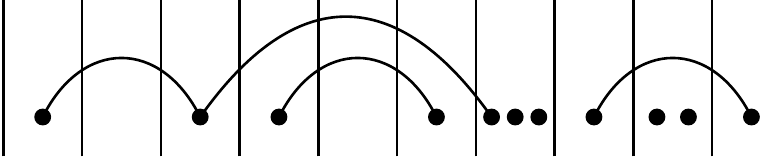}}
\]
Assume now two consecutive entries on one side are
equivalent. The above restrictions imply that
\[
\cl{\includegraphics[scale=0.9]{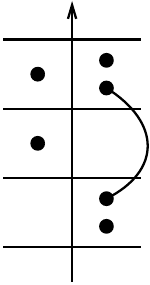}}
\]
the first one is the last of some row (on some side), and the next one
is the first of the second next row (on the same side), with the row
in between being empty (on the same side).

If in particular three consecutive entries are equivalent,
then the middle one is single on its side and row.

Consider the following transformation on the equivalence arcs
(for a moment disregarding the row structure) on both sides.
\begin{eqn}\label{(5)}
\lower15pt\hbox{\includegraphics[scale=0.9]{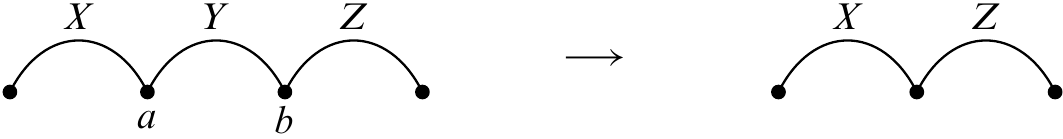}}
\end{eqn}
This transformation eliminates a middle one in a series of
$4$ consecutive equivalent entries (that is, below the arcs
the space is assumed to be empty). We will call an arc like
the middle arc $Y$ on the left side and the beads
$a,b$ it connects \emph{shrinkable}. Arcs and beads not representable in
this form are called \emph{non-shrinkable}. We will apply
\eqref{(5)}, but before this let us count how many times
we can do so.

\begin{sublemma}
At most $\!\myfrac{3}{8}$ of the beads in the table
(counted altogether on \emph{both} sides of the table)
are removable by the transformation \eqref{(5)}.
\end{sublemma}

(Note that there are more shrinkable beads than removable ones,
as on the right of \eqref{(5)} a shrinkable bead remains.)

\proof Consider a shrinkable arc $Y$ in \eqref{(5)}. Such an arc
identifies entries $a$ and $b$ two rows apart, so that 
$a$ and $b$ are single on their row and side, and the row between
them is empty (on that side).

Let $A,B$ be the slots of the table in the rows of $a$ and $b$,
respectively, but on the opposite side to $a$ and $b$. Let $C$ be
the slot in between $A$ and $B$.
\begin{eqn}\label{(11)}
\lower50pt\hbox{\includegraphics[scale=0.9]{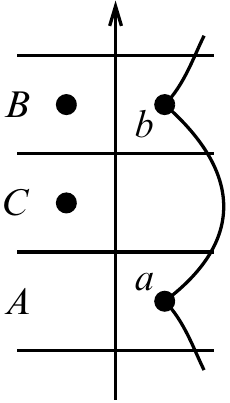}}
\end{eqn}
Then by the restrictions there are two alternatives:

\renewcommand\theenumi{\alph{enumi}}
\begin{enumerate}
\item\label{item&a} $C$ contains at least two entries or
\item\label{item&b} $C$ contains one entry, and at least
  one of $A$ and $B$ contains at least one entry.
\end{enumerate}

We would like to count how many entries/beads in $A,B,C$ for all 
shrinkable arcs in \eqref{(5)} we find this way, such that the beads
themselves are non-shrinkable. To ensure the correct counting,
we must take for each occurrence of \eqref{(5)} each such bead in
$A$, $B$ or $C$ with weight $1/k$, where $k$ is the number of different
fragments \eqref{(5)} for which the bead occurs in \eqref{(11)}.
Our aim is to group non-shrinkable beads and shrinkable arcs
so that the weighted bead count is at least $5/3$ per shrinkable arc.

It is easy to see that in both of the above cases any of the entries
in $A$, $B$ or $C$ are non-shrinkable,
and that the entries on $C$ are counted once, while those on $A$ or
$B$ are counted (once or) at most twice for different fragments \eqref
{(11)}. This means that the contribution of \eqref{(5)} to the
above weighted count of $A,B,C$ beads is \emph{at least two}, unless
only one bead is on $A$ or $B$, and this bead is counted twice (ie
$k=2$).

Assume thus the \emph{only} entry on $A$ or $B$ is counted twice.
We have then the following table fragment:
\begin{eqn}\label{65}
\lower50pt\hbox{\includegraphics[scale=0.9]{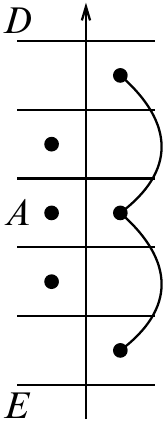}}
\end{eqn}
Since $a$ and $b$ are shrinkable, there must be one further
row on top (and bottom), whose slot on the opposite side to $a$ and $b$
we call $D$ (and $E$), and $D$ (resp.\ $E$) must contain
$\ge 2$ beads. (The slot in the row of $D,E$ on the side of
$a$ and $b$ is empty because there is at least a terminal arc
as $X,Z$ in \eqref{(5)}.) These elements may be counted
for another shrinkable arc, or not. If there is
such an arc, it is on the same side as $a$ and $b$,
and part of a longer sequence of consecutive
shrinkable arcs (ie $X$ or $Z$ in \eqref{(5)} is also shrinkable).
\[
\cl{\includegraphics[scale=0.9]{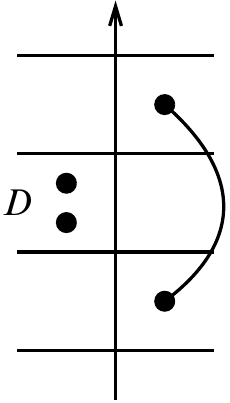}}
\]
We have thus two or three shrinkable arcs, and have found to them
at least 5 entries. It is easy to see that all of these
entries are non-shrinkable
(because shrinkable arcs connect only single entries on their row,
two rows apart, with either neighboring rows being empty
on their side of the table). Also any of these entries is
not counted again for another shrinkable arc in the above
grouping, \emph{unless} we have:
\[
\cl{\includegraphics[scale=0.9]{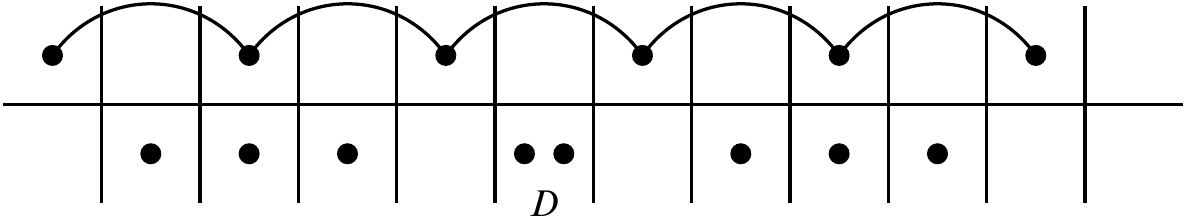}}
\]
\newpage
(The slot of $D$ may have more than 2 beads.) Then we are in the
situation of \eqref{65}, just with 3 more shrinkable arcs and
$\ge 5$ more non-shrinkable beads. We could then argue inductively
by the number of consecutive shrinkable arcs.

In summary we found either for a shrinkable arc
two non-shrinkable entries (counted with correct weight),
or to 2 or 3 shrinkable arcs, 5 non-shrinkable entries
(all with weight 1). Thus if there are $l$ arcs/entries
in the table removable by \eqref{(5)},
at least $\myfrac{5}{3}\,l$ other entries are not removable. \endproof

Now discard distinction between the left and right side, and
the subdivision into rows, and consider
the pictures for each side of the table separately.
In the remaining picture after the transformation
\eqref{(5)} we have no 4 consecutive equivalent beads 
on the same side. 

Consider the equivalence class of the first element. Let it
have $t$ elements.
\begin{eqn}\label{(6)}
\lower 30pt\hbox{\includegraphics[scale=0.9]{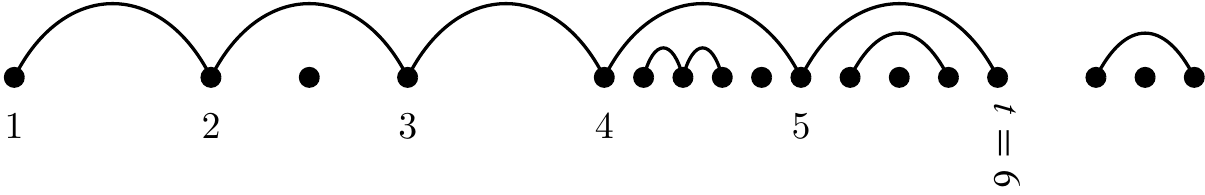}}
\end{eqn}
Removing this class, and cutting along the $t$ points,
we have $t$ pieces satisfying the following conditions:
\begin{itemize}
\item no portion contains 4 consecutive equivalent
elements,
\item from the first $t-1$ portions, at least every third
is non-empty.
\end{itemize}

Let $nr(n)$ be the minimal number of equivalence classes
of $n$ elements under this equivalence relation. We claim

\begin{sublemma}
\[
nr(k)\ge\left\{\begin{array}{cc}\frac{k}{6}+\frac{1}{2} & k>0 \\
0 & k=0 \end{array}\right\}\,.
\]
\end{sublemma}

\proof We proceed by induction on $k$.
We have $nr(0)=0$ and $nr(1)=nr(2)=nr(3)=1$. If $k>3$,
we find from the above picture \eqref{(6)} the recursion
\begin{eqnarray}
\nonumber
nr(k) & \ge & 1+\min\, \left\{
\es\sum_{i=1}^t\,nr(l_i)\es:\es
\vcbox{\shortstack[c]{$t>0,\,\ l_1+\dots+l_t=k-t,\,\ l_i\ge
0,$\,\\ at least $\br{\frac{t-1}{3}}$
of the $l_i>0$}}
\es\right\} \\[2mm]
&\nonumber \ge &
1+\frac{k-t}{6}+\br{\frac{t-1}{3}}\cdot
\frac{1}{2}\es\ge\es
1+\frac{k-t}{6}+\left(\frac{t}{3}-1\right)\cdot
\frac{1}{2}\es=\es\frac{k}{6}+\frac{1}{2}\,.\rlap{\hspace{-3.3mm\qed}}
\end{eqnarray}

We have $r:=\len\gm'-1$ rows in the table. %For each row,
As remarked, any third row in the table has at least two entries.
Thus the total number of entries is at least
\[
\frac{4r-2}{3}\,=\,\frac{4(\len\gm'-1)-2}{3}\,=\,\frac{4}{3}\len\gm'-
2\,.
\]
At most $\myfrac{3}{8}$ of these entries
were discarded under the identification
\eqref{(5)}. By applying the last sublemma, we see that
from the remaining at least $\myfrac{5}{8}\left(
\frac{4}{3}\len\gm'-2\right)$ entries, at least 
$\myfrac{1}{6}$ plus one more remained after identifying regions
(note that none of the sides of the table was empty, even after 
\eqref{(5)}).

Thus there are at least
\[
\frac{1}{6}\cdot\frac{5}{8}\left(
\frac{4}{3}\len\gm'-2\right)+1\,=\,
\frac{5}{48}\left(\frac{4}{3}\len\gm'-2\right)+1
\]
other regions not passed by $\gm'$, which have neighbors passed
intermediately by $\gm'$. We had in $D'$ two more regions, 
$D'_{1,2}$, not passed by $\gm'$ with no neighbored region
passed intermediately by $\gm'$. Thus there are at least
\[
\BR{\frac{5}{48}\left(\frac{4}{3}\len\gm'-2\right)+3}\,=\,
\BR{\frac{5}{36}\len\gm'+\frac{67}{24}}\,\ge\,\frac{5\,\len\gm'+101}{36}
\]
regions not passed by $\gm'$. Then $\gm'$ passes $\len\gm'+1$ other
regions, so
\[
r(D')\,\ge\,\frac{41}{36}(\len\gm'+1)\,+\frac{8}{3}\,,
\]
which leads to the first maximum alternative of the inequality
we claimed in \eqref{mxz}. The Curve length lemma is now proved. \endproof

\section{Genus decreasing bound}

Once we have a bound on the number of regions a curve $\gm$
of the above specified type passes, we use this bound
to obtain a bound on the decrease of genus of a positive
knot diagram we can achieve by replacing $\gm$ by $\gm_1$.
In addition to the wave moves, we need the move:
\begin{eqn}\label{99}
\lower 20pt\hbox{\includegraphics[scale=0.9]{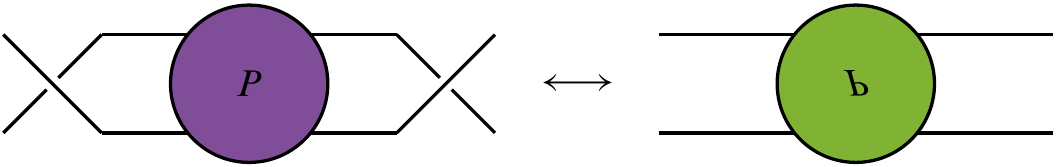}}
\end{eqn}
Since for the empty tangle $P$ this is a clasp resolution, we call
\eqref{99} a generalized clasp reduction move. We can consider a
clasp resolution as a special case, included in \eqref{99}.
\newpage
We assume that if $P$ is empty, and the r.h.s.\ of \eqref{99} has
nugatory crossings, that we remove all these crossings
(ie reduce the diagram).

\begin{theorem}\label{thgdec}
Let $D$ be a positive diagram of a knot $K$ with $g=g(D)=g(K)>1$.
\vspace{2pt}

Then by applying switches of crossings in $D$, and subsequently a wave
move (shrinking a bridge/tunnel to length 1), or a move \eqref{99}, we
can obtain a positive diagram $D'$ of a positive knot $K'$ such that
\emph{either} $g(K')=g(K)$ and $c(D')<c(D)$, \emph{or} $g(K')<g(K)$ and
\begin{eqn}\label{geq}
g(K')\,\ge\,\frac{6}{13}+ \frac{41}{221}g(K)\,.
\end{eqn}
\end{theorem}
\vspace{2pt}

For the proof we need to quote one more previous result, about
the maximal crossing number $c_g$ of generators (see \fullref{df1}). In \cite{gen1}, a rather rough estimate on
the number $d_g$ was given, which was later improved in
\cite{STV} to $d_g\le 6g-3$. Then in \cite{SV} we showed
that this inequality is sharp. Since $c_g\le 2d_g$, we have
$c_g\le 12g-6$. Later, in \cite{SV}, we showed by explicit
examples that $c_g\ge 10g-7$, and remarked that at least $c(D)
\le 10g-6$ if $D$ is a special alternating generator. The work was
completed in \cite{adeq}, where, using Hirasawa's algorithm,
the maximal generator crossing number was determined (also for links).
There the value $10g-7$ was found exact, also for arbitrary generators.
\vspace{2pt}

\begin{theorem}\label{thc_g}{\rm\cite{adeq}}\qua
Assume $g\ge 2$. Then $c_g=10g-7$. Moreover, generators of
genus $g$ with the maximal number of crossings are always
special alternating.
\end{theorem}
\vspace{2pt}

\proof[Proof of \fullref{thgdec}]
If $D$ has a bigon region (clasp), we are easily done, since
$g(K')\ge g(K)-1$. So assume $D$ has no such region.
\vspace{2pt}

It is easy to see that it suffices also to prove the theorem in case
$D$ is prime; the composite case follows easily. For connected sum
factors of genus $1$ one can use the description in \cite{gen1};
since all diagrams have a clasp, the genus of $D$ would decrease
at most by one. 
\vspace{2pt}

Choose a minimal admissible
curve $\gm$ in $D$ of length (number of edges in $D\sm\gm$
it intersects) $\len\gm>1$, such that start and end region are
neighbored (in $D\sm\gm$). In the proof of the Curve length lemma we
argued that an admissible curve always exists. Among these curves $\gm$
take one of minimal length. Then this curve is one of the type
considered in that lemma.
\vspace{2pt}

So now we can apply the Curve length lemma. Assume first the first
alternative in the maximum of \eqref{mxz} holds. Since $\gm$ passes
$\len\gm+1$ regions (including first and last), we can write the
estimate as
\[
\len\gm\,\le\,\frac{36r(D')-137}{41}\,,
\]
where $r(D')$ is the number of regions of $D'$, and $D'$ is
the diagram obtained by repacing $\gm$ in $D$ by a
one-crossing curve $\gm_1$. 

We have
\[
r(D')\,=\,2+c(D')\,=\,2+c(D)-\len\gm+1\,=\,3+c(D)-\len\gm\,.
\]
Thus
\[
\len\gm\,\le\,\frac{36(3+c(D)-\len\gm)-137}{41}\,=\,
\frac{36(c(D)-\len\gm)-29}{41}\,,
\]
and
\begin{eqn}\label{(one)}
\len\gm\,\le\,\frac{36c(D)-29}{77}\,.
\end{eqn}
We can achieve $D'$ to be positive by properly switching
crossings of $\gm$ to become under- or over-crossings,
dependingly on how the one crossing of $\gm_1$ is to be
switched to be positive. Since $\gm\cup\gm_1$
gives a closed curve, which can be perturbed to
be transversal to $D$ without altering the
edges $\gm$ and $\gm_1$ intersect, a linking number argument
shows that the number of crossings on $\gm$ to
be switched is 
\[
\frac{\len\gm}{2}\mbox{\ \ or\ \ }\frac{\len\gm\pm 1}{2}\,.
\]
Assume the largest possible value, since we want to have
only an upper bound on the number of such crossings. Thus,
by switching at most {\small$\ds\frac{\len\gm+1}{2}$}
crossings in $D$, we can simplify it to a positive diagram $D'$.

We would like to show now that $D'$ is the diagram whose existence
was asserted in the theorem.

First, clearly $c(D')<c(D)$. To see $g(D')\le g(D)$, use
that $g(D)=g(K)$ and $g(D')=g(K')$. Then consider the diagram
$D''$ obtained from $D$ after the crossing switches on $\gm'$,
but before the wave move taking it into $D'$ (replacing
$\gm$ by $\gm_1$). This is certainly a diagram of $K'$, so
that $g(K')\le g(D'')=g(D)$.

On the other hand, by applying Bennequin's inequality \eqref{BQ}
on the diagram $D''$ of $K'$, we have
\[
g(K')\,=\,g(D')\,\ge\,g(D)-\frac{\len\gm+1}{2}\,.
\]
Then, applying \eqref{(one)}, we obtain
\begin{eqn}\label{(two)}
g(D')\,\ge\,g(D)-\frac{18c(D)}{77}-\frac{24}{77}\,.
\end{eqn}
There is now the possibility to consider the other alternative
of the maximum in \eqref{mxz} holds. In this case $\len\gm\le 3$, so
we can switch ($1$ or) $2$ crossings, and have a transformation
into a diagram $D'$ with $c(D')<c(D)$. By the same argument using
\eqref{BQ} as above, $g(D)\ge g(D')\ge g(D)-2$. The second inequality
fits into \eqref{(two)}, except if $c(D)\le 7$. But such $D$ always
has a clasp (since it is a knot diagram; by direct verification),
which we excluded. So we can ignore the second maximum alternative
in \eqref{mxz}, and use \eqref{(two)}.

The next step is to apply generator estimates. As a preparation,
we need to establish two properties of $D'$.

\begin{sublemma}
$D'$ has no nugatory crossing.
\end{sublemma}

\proof Assume that $c$ were such a crossing.
\[
\cl{\includegraphics[scale=0.9]{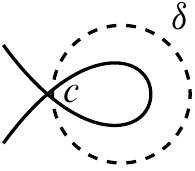}}
\]
Then there is a nugatory curve $\dl$. \fullref{Lm4} implies
that $\gm'$ cannot pass twice the region $X$ of $D'$ that contains
$\dl$ (as in the right picture below). So $\gm'$ remains (as in the
left picture) within $X$, and can be isotoped off $\dl$ (within $X$).
Then $c$ becomes nugatory in $D$, too, a contradiction to its
assumed primeness.
$$
\lower5pt\hbox{\includegraphics[scale=0.9]{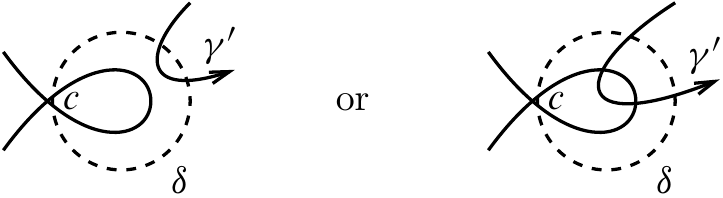}}
\proved$$

\begin{sublemma}
$D'$ admits at most one reducing $\bm$ move.
\end{sublemma}

\proof Assume there are $3$ $\sim$--equivalent
crossings in $D'$. Consider first the case that
the crossing of $\gm_1$ is not among them.
Let $A$ and $B$ be the (common) non-Seifert circle
regions at these crossings. Let $\ap$, $\bt$, $\dl$
be the curves connecting some fixed point in $A$ with
some fixed point in $B$ via one of these crossings.
Then the join of any two of these curves gives a
closed curve intersecting $D'$ transversely in
$2$ crossings:
\[
\lower45pt\hbox{\includegraphics[scale=0.9]{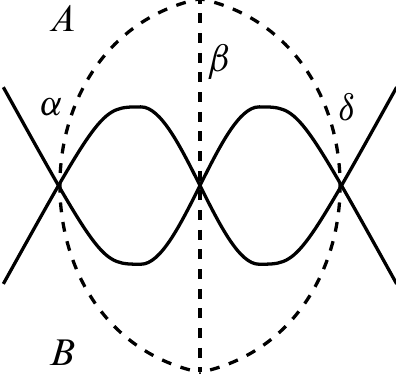}}
\qquad\text{\mbox{in $D'$.}}
\]
Now we delete $\gm_1$, which does not affect the three
crossings and curves, and try to reinstall $\gm$ so
as to obtain $D$.

Since $\gm'$ passes every region of $D'$ at most once,
if $\gm'$ passes through $\ap\cup\bt$ without being homotopable
off it,
\[
\cl{\includegraphics[scale=0.9]{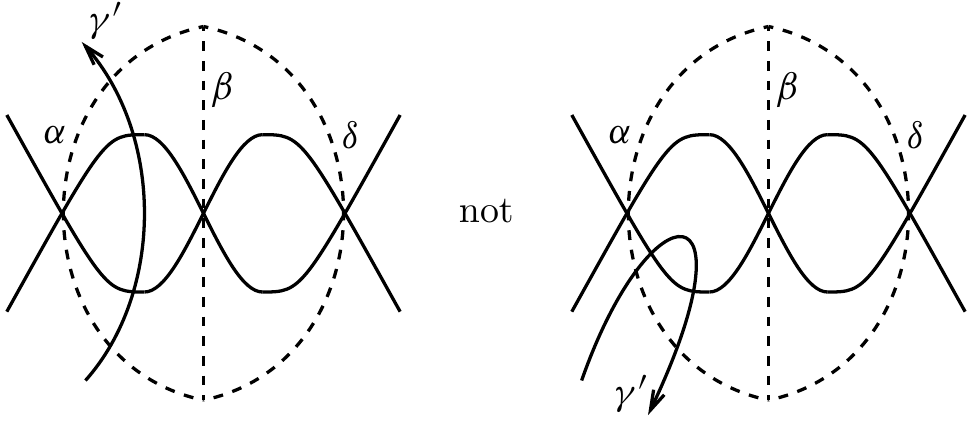}}
\]
then it cannot pass through $\bt\cup\dl$. Thus up to flype, we have
in $D$ a clasp, and can proceed by one crossing change and
\eqref{99}. What appears as a clasp in the above picture may in
fact be a non-trivial tangle. Thus there is also the option that
$\gm'$ starts and ends within one of these parts:
\[
\cl{\includegraphics[scale=0.9]{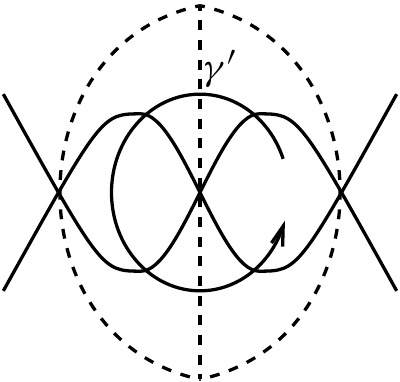}}
\]
In this case it is $\ap\cup\dl$ which is not affected (up to homotopy),
and then as above \eqref{99} applies.
\newpage

However, in the second case
when one of the three crossings is the one on $\gm_1$,
one can have a curve $\gm'$ so that $\gm$ makes in $D$ all the
$\sim$--equivalent crossings inequivalent:
\begin{eqn}\label{xx}
\cl{\includegraphics[scale=0.9]{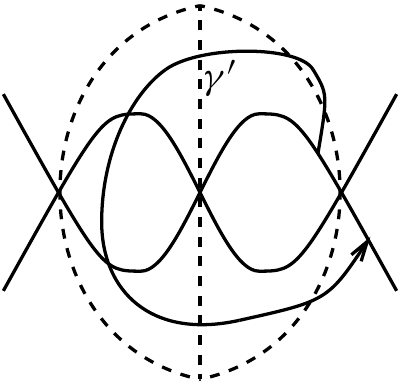}}
\end{eqn}
Since there is only one crossing on $\gm_1$, the situation
in \eqref{xx} can happen in at most one $\sim$--equivalence class
of $D'$, and if there are 4 or more crossings in that class,
again \eqref{99} can be used. \endproof

Now use our previous work on $c_g$, the maximal number of crossings
of a genus $g$ reduced $\bm$ irreducible knot diagram
(generator) of genus $g$. We have from \fullref{thc_g}
that $c_g\le 10g-7$. Therefore, with the previous two sublemmas,
we have $c(D')\le 10g(D')-5$, and
\begin{eqnarray*}
g(D') & \ge & \frac{c(D')+5}{10} \, = \, \frac{c(D)-\len\gm+6}{10} \\
& \stackrel{\eqref{(one)}}{\ge} &
\frac{\ds c(D)-\frac{36c(D)-29}{\phantom{,}77}+6}{10} \, = \, 
\frac{41}{770}c(D)+\frac{491}{770}\,.
\end{eqnarray*}
This, together with \eqref{(two)}, gives
\begin{eqn}\label{(star)}
g(D')\,\ge\,\max\,\left\{\,
g(D)-\frac{18c(D)+24}{77},\,\frac{41}{770}c(D)+\frac{491}{770}
\right\}
\,.
\end{eqn}
To determine the worst possible value of this estimate for
fixed $g(D)$, one has to equate both alternatives:
\begin{eqnarray*}
g(D) & = & \frac{18c(D)+24}{77} + \frac{41}{770}c(D)+\frac{491}{770}
\, = \, \frac{221}{770}c(D)+\frac{731}{770}\,,\qquad\mbox{and so} \\
c(D) & = & \frac{770}{221}g(D)-\frac{731}{221}\,.
\end{eqnarray*}
Then the r.h.s.\ of \eqref{(star)} evaluates to
\[%\begin{eqnarray*}
g(D') \, \ge \, \frac{41}{770}c(D)+\frac{491}{770} \, = \,
\frac{41}{770}\,\left(\frac{770}{221}g(D)-\frac{731}{221}\right)+
\frac{491}{770} 
\, = \, \frac{6}{13}+\frac{41}{221}g(D)\,,
\] %end{eqnarray*}
which gives the asserted estimate \eqref{geq}. \qed

\section{Applications}

\subsection[Proof and extension of \ref{thalg}]{Proof and extension of \fullref{thalg}}

We are now done with most of the work concerning \fullref{thalg}. Define for $n>1$ the numbers $g_n$ by
\begin{eqn}\label{dgn}
g_n\,=\,\ds\br{\frac{221n-323}{41}}\,,
\end{eqn}
which we will use throughout this section.
Note that $n\le g_n$ when $n\ge 2$. As a consequence of
\fullref{thgdec}, we obtain:

\begin{prop}\label{pp42}
For $n>1$, let $X_n$ be the set of positive generating
diagrams of genus $n,\dots,g_n$.
Let $D$ be a positive knot diagram of genus at least $n$.
Then there is a sequence of diagrams $D_0,D_0',D_1,D_1',\dots,
D_{k-1}',D_k$ with $D_0=D$ and $D_k\in X_n$, satisfying the
following properties:
\begin{enumerate}
\item $D_i$ is positive and $D_i'$ differs from $D_i$ by
crossing changes (ie, $D_i$ and $D_i'$ have the same shadow),
\item $D_{i+1}$ is obtained from $D_i'$ either by a move \eqref{99},
or by rerouting a bridge/tunnel to length 1\,.
\end{enumerate}
\end{prop}

\proof 
% One argues by induction on $n$, for given $n$ over $g(D)$,
% and for fixed $g(D)$, by induction on $c(D)$. 
% 
Let $D$ be a positive diagram of genus $\ge n$. If $g(D)\le g_n$, then
we need to reduce $D$ only to the generating diagram in whose series
it lies, and for this the move \eqref{99} (after the proper
crossing changes) is enough. So assume $g(D)>g_n$.
% 
% If $D=D_1\#\dots\#D_l$ is composite, argue for its prime factors $D_i$
% separately, using induction over $n$. Choose numbers $n_i$ with
% $0\le n_i\le g(D_i)$ and $\sum n_i=n$. If $n_i>1$, one can reduce
% by induction $D_i$ to a positive diagram $D_{i,k_i}\in X_{n_i}$. If $n_i
% \le 1$, one can reduce $D_i$ either to an unknot or a (positive) trefoil
% diagram $D_{i,k_i}$. Using the superadditivity of $g_n$ in $n$ when $n\ge
% 2$ and adjusting a suitable proportion of the $D_i$ with $n_i\le 1$
% to change to unknots and trefoils, we can achieve that the genus of
% the resulting diagram $D_{1,k_1}\#\dots\#D_{l,k_l}$ lies between $n$ and $g_n$.
% 
% So assume $D$ is prime.
% Because $g(D)>g_n$ and
Since $g(D)\in\bZ$, we have then
% $g(D)\ge\mbox{\small$\ds\frac{221n-322}{41}$}$, and
\[
g(D)> \frac{221n-323}{41}\,=\,\frac{221(n-1)-102}{41}\,.
\]
Then by \fullref{thgdec} we can switch crossings in $D=D_0$
to a diagram $D_0'$, which we can transform into a positive diagram
$D_1$, such that $c(D_1)<c(D)$, $g(D_1)\le g(D)$, and
\[
g(D_1)\,\ge\,\frac{6}{13}+\frac{41}{221}g(D)\,>\,n-1\,.
\]
So $g(D_1)\ge n$, and we can argue with $D_1$ by induction on $g(D)$,
and for fixed $g(D)$, inductively over $c(D)$. \endproof

\begin{corr}\label{cgn}
For $n>1$ and $g_n$ as in \eqref{dgn} we have
\begin{align}\label{sig=}
\min\,\left\{\,\sg(K)\,:\,\mbox{$K$ positive},\ n\le g(K)\le
g_n\ \right\}\es&=\\
\min\,\{\,\sg(K)\,:\,&\mbox{$K$ positive},\ n\le g(K)\ \}\,.\nonumber
\end{align}
In particular, assume that there is no positive knot of
signature $\le \sg$ and genus $n,\,\dots,\,g_n$\,. Then there is
no positive knot of signature $\sg$ and genus $\ge n$ at all.
\end{corr}

\proof This follows straightforwardly from \fullref{pp42}
and the property \eqref{2a} of $\sg$. \endproof

\proof[Proof of \fullref{thalg}]
Again this is an immediate consequence of \fullref{pp42},
as in \fullref{cgn}: let $C_n$ be the set of knots with 
positive generating diagrams of genus $n$ to $g_n$. \endproof

\begin{rem}
The estimates may be improvable by 
using more consequently the integrality of $c(D)$ and $g(D)$,
which is not always guaranteed in the above calculations.
This improvement, however, will affect only the absolute term,
and expectedly not in a significant way, so that we preferred
to largely waive on incorporating this additional effort
into our (anyway technical enough) proofs.
\end{rem}

More generally, we have:

\begin{theorem}\label{th5.3}\renewcommand\theenumi{\arabic{enumi}}\
\begin{enumerate}
\item\label{item__1}
Let $v$ be an invariant with the following property:
if $D$ is a positive diagram and $D'$ is obtained by (some
non-zero number of) crossing switches from it,
then $v(D')<v(D)$. Then $v$ has an increasing
lower bound on positive knots in terms of the
genus of the knot.
\item\label{item__2}
If $v(D')\le v(D)$, then one can algorithmically partially decide
whether any given value of $v$ is attained on positive knots of only
finitely many genera.
\end{enumerate}
\end{theorem}

Examples of invariants satisfying $v(D')<v(D)$ are the
properly scaled Vassiliev invariants $v_2$ and $v_3$ 
of degree 2 and 3. 
In fact, we showed in \cite{pos} independently that for
either of $v_2$ and $v_3$ both the premise and conclusion of
part \ref{item__1} of the above theorem hold.
(For $v_2$ we must exclude $D'$
being the mirror image of $D$, but the argument for
\fullref{th5.3} clearly still works under
this small restriction.)
It is interesting to reveal that there is in fact
a certain causality between both results, which is
quite non-evident from the approach in \cite{pos}.

Note also that the recent signature-type concordance invariants
of Ozsv\'ath and Szab\'o \cite{OS} and Rasmussen \cite{Rasmussen}
satisfy $v(D')\le v(D)$, though the conclusion we obtain is known.
More interestingly, one can apply part \ref{item__2} of the
theorem to the Tristram--Levine signatures, thereby extending our
treatment of $\sg$. 

\subsection{Partial orders of knots}

To conclude with, we mention a relation to the partial orders of
Taniyama and Cochran--Gompf.
In \cite{Taniyama}, Taniyama defined a partial order of knots
by $K_1\succeq K_2$ if the set of shadows of diagrams of $K_1$ is a
(not necessarily proper) superset of the set of shadows of
diagrams of $K_2$. 

The standard fact
that any diagram is unknottable by crossing changes means
that the unknot is the maximal element in this partial order.
Taniyama showed that the trefoil dominates any other
knot, and that $5_1$ dominates any knot except connected sums
of $(p,q,r)$--pretzel knots, $p,q,r$ odd. The first result
gives an easy proof that positive knots have positive
signature, and the second result amplifies this statement
by showing that $\sg\ge 4$ if $g\ge 2$. What our arguments show
can be thought of as a generalization of Taniyama's two results.

\begin{defi}\label{xdf}
Call a set $\cK_1$ of knots \emph{dominating} a set $\cK_2$, in notation
$\cK_1\succeq\cK_2$, if the set of shadows of all diagrams of knots
in $\cK_1$ is a (not necessarily proper) superset of the set of
shadows of diagrams of all knots in $\cK_2$.
A set $\cK$ of knots is \emph{finitely dominated}, if there is
a finite subset $\cK'\subset \cK$ dominating $\cK$.
\end{defi}

\begin{theorem}\label{thT}
For all $n\ge 1$, the set $G_n$ of knots of canonical genus $\ge n$
is finitely dominated. There is a finite subset $C_n\subset G_n$
of positive knots, such that $C_n\succeq G_n$.
\end{theorem}
% \vspace{4mm}

\begin{rem}\label{rzq}
It is clear that for two knots $K_{1}\succeq K_2$ implies
$c(K_1)\le c(K_2)$. Thus any chain of `$\succeq$' has a maximal
element. This means also that the subset of maximal elements is
a dominating subset. It is not clear, though, that conversely a
dominating subset must contain (all, or even any) maximal elements.
\fullref{xdf} does \emph{not} imply that for each $K\in\cK_2$
there is a $K'\in\cK_1$ with $K'\succeq K$, because the $K'$ we
find for different diagram( shadow)s of $K$ may not be the same.
Therefore, we do not know if the set $G_n$ in \fullref{thT}
has only finitely many maximal elements. (See, however, in contrast
\fullref{thY}.)
\end{rem}

\proof[Proof of \fullref{thT}]
Consider $n\ge 3$, since $C_1=\{3_1\}$ and
$C_2=\{5_1,\,3_1\#3_1\}$ are Taniyama's results.

With $g_n$ defined as in \eqref{dgn}, let the
set $C_n$ consist of the positive knots with positive generating
diagrams of genus $n,\dots,g_n$, and write $X_n$ for
the set of these diagrams. Then apply an induction argument
(similar to the one in \cite[remark 6.2]{apos}).

By \fullref{pp42}, for any diagram $D$
of genus $\ge n$, there is a sequence of diagrams
$D_0,\dots,D_k=D$ (we reindexed the subscripts here) such that
$D_0\in X_n$, and $D_i$ and $D_{i+1}$ differ by crossing
changes, followed either by a generalized clasp resolution
\eqref{99} or a wave-move of a bridge/tunnel to one of length 1.
The problem is how to switch crossings in $D=D_k$, so that one
can perform all these moves without switching later crossings
in between the moves.

To see how to do this, we use induction on $k$. The case $k=0$ is
clear. If $k>0$, we know by induction that we can switch crossings
in $D_{k-1}$ so that it reduces to $D_0$. The case that
$D_{k-1}$ is obtained from $D_k=D$ by the move \eqref{99},
it is clear how to choose the crossing switch of $D_k$. (Note that
the crossings in the tangle $P$ in \eqref{99} are flipped around.)

So consider the case that $D_{k-1}=D'$ is obtained from $D_k=D$ by
shrinking a bridge/tun\-nel. Then we have an arc $C$ in $D$ which can be
shrunk by a wave-move to an arc $C'$ in $D'$ of a single crossing $p$.
Since one can crossing-switch $C$ prior to the wave-move to pass
above or below the rest of the diagram, one can adjust the sign
of the crossing $p$ it collapses to. This crossing $p$ in $D'=D_
{k-1}$ will be switched at most once in the simplification from
$D_{k-1}$ to $D_0$. If it is switched (resp.\ not switched) so that the
strand of $C'$ becomes (resp.\ remains) an over/undercrossing in $D'$,
then switch $C$ so that it becomes a bridge/tunnel in $D$.
\endproof

\begin{rem}
It is clear that one can choose the set $C_n$ to consist of
alternating knots, instead of positive ones.
\end{rem}

A different, but related partial order of knots was introduced by 
Cochran and Gompf \cite{CochranGompf}. In their sense, $K_1\ge K_2$
if $K_1$ is concordant to $K_2$ inside a 4--manifold with positive
intersection form. This occurs for example if $K_2$ is obtained
from $K_1$ by changing a positive crossing to a negative one. The
referee pointed us to make a comment in that context, which we finish
with. Clearly, as for Taniyama's partial order, our work directly
connects also to the one of Cochran--Gompf. From \fullref{pp42} we obtain immediately:

\begin{theorem}\label{thY}
For each $n\ge 1$, the set of positive knots of genus at least $n$ has
a finite number of minimal elements in Cochran--Gompf's partial order.
\end{theorem}

The case $n=1$ is not formally included in the proposition, but is
known.
Cochran and Gompf showed that a non-trivial positive knot dominates
the right-handed trefoil $!3_1$, which is also a simple implication
of Taniyama's related theorem. His second theorem similarly
shows that any positive knot of genus $>1$ dominates in `$\ge$'
one of $!5_1$ or $!3_1\#!3_1$. In general it is easy to see that
if $K_1$ is positive, then $K_2\succeq K_1$ implies $K_1\ge K_2$.
The converse implication is false even if we assume that $K_2$
is also positive; $K_1=!5_1$ and $K_2=!5_2$ give an example.
(Slightly fancier pairs of a 12 crossing knot $K_1$ and a 13
crossing knot $K_2$ show that `$\ge$', even for positive knots,
also fails to respect the crossing number~-- in contrast to
`$\succeq$'; see \fullref{rzq}.)
Still a more detailed study of the relationship between `$\ge$'
and `$\succeq$', which seems not to have been undertaken
so far, may be worthwhile.

\noindent{\bf Acknowledgement}\qua This paper emerged from work I
carried out over a long period at several places. Some part of
this work was written during my stay at the
University of Toronto in spring 2002, supported by a grant of
Deutsche Forschungsgemeinschaft (DFG). The paper was completed
at the Graduate School of Mathematical Sciences, University of Tokyo.
Its hospitality, the support by Postdoc grant P04300 of Japan
Society for the Promotion of Science (JSPS) and of my host Professor
T~Kohno is also acknowledged.
The author is supported by 21st Century COE Program.

\bibliographystyle{gtart}
\bibliography{link}

\begin{thebibliography}{}
\providecommand\bibmarginpar{\leavevmode\marginpar}
\def\urlstyle#1{{\tt #1}}

\bibitem{ACampo}
\textbf{N A'Campo}, \href{http://www.numdam.org/item?id=PMIHES_1998__88__151_0}
  {\emph{Generic immersions of curves, knots, monodromy and {G}ordian number}},
  Inst. Hautes \'Etudes Sci. Publ. Math.  (1998) 151--169 (1999)
  \xox{MR}{1733329}

\bibitem{Bennequin}
\textbf{D Bennequin}, \emph{Entrelacements et \'equations de {P}faff}, from:
  ``Third Schnepfenried geometry conference, Vol. 1 (Schnepfenried, 1982)'',
  Ast\'erisque 107, Soc. Math. France, Paris (1983)  87--161 \xox{MR}{753131}

\bibitem{BoiWeb}
\textbf{M Boileau}, \textbf{C Weber}, \emph{Le probl\`eme de {J}. {M}ilnor sur
  le nombre gordien des n\oe uds alg\'ebriques}, Enseign. Math. $(2)$ 30 (1984)
  173--222 \xox{MR}{767901}

\bibitem{Brittenham}
\textbf{M Brittenham}, \emph{Bounding canonical genus bounds volume},\ \
  preprint (1998)
\ Available at \setbox0\hbox{\makeatletter\@url
{http://www.math.unl.edu/~mbritten/personal/pprdescr.html}}
\href{http://www.math.unl.edu/~mbritten/personal/pprdescr.html}
{\unhbox0}

\bibitem{Busk}
\textbf{J\,M van Buskirk}, \emph{Positive knots have positive {C}onway
  polynomials}, from: ``Knot theory and manifolds (Vancouver, 1983)'', Lecture
  Notes in Math. 1144, Springer, Berlin (1985)  146--159 \xox{MR}{823288}

\bibitem{Cochran}
\textbf{T\,D Cochran}, \href{http://dx.doi.org/10.2140/agt.2004.4.347}
  {\emph{Noncommutative knot theory}}, Algebr. Geom. Topol. 4 (2004) 347--398
  \xox{MR}{2077670}

\bibitem{CochranGompf}
\textbf{T\,D Cochran}, \textbf{R\,E Gompf},
  \href{http://dx.doi.org/10.1016/0040-9383(88)90028-6} {\emph{Applications of
  {D}onaldson's theorems to classical knot concordance, homology {$3$}-spheres
  and property {$P$}}}, Topology 27 (1988) 495--512 \xox{MR}{976591}

\bibitem{Conway}
\textbf{J\,H Conway}, \emph{An enumeration of knots and links, and some of
  their algebraic properties}, from: ``Computational Problems in Abstract
  Algebra (Proc. Conf., Oxford, 1967)'', Pergamon, Oxford (1970)  329--358
  \xox{MR}{0258014}

\bibitem{Cromwell2}
\textbf{P\,R Cromwell}, \emph{Homogeneous links}, J. London Math. Soc. $(2)$ 39
  (1989) 535--552 \xox{MR}{1002465}

\bibitem{Cromwell}
\textbf{P\,R Cromwell}, \emph{Positive braids are visually prime}, Proc. London
  Math. Soc. $(3)$ 67 (1993) 384--424 \xox{MR}{1226607}

\bibitem{Crowell}
\textbf{R Crowell},
  \href{http://links.jstor.org/sici?sici=0003-486X(195903)2:69:2%3C258:GOALT%3%
E2.0.CO%3B2-B} {\emph{Genus of alternating link types}}, Ann. of Math. $(2)$ 69
  (1959) 258--275 \xox{MR}{0099665}

\bibitem{FintushelStern}
\textbf{R Fintushel}, \textbf{R\,J Stern},
  \href{http://links.jstor.org/sici?sici=0003-486X(198509)2:122:2%3C335:PO%3E2%
.0.CO%3B2-H} {\emph{Pseudofree orbifolds}}, Ann. of Math. $(2)$ 122 (1985)
  335--364 \xox{MR}{808222}

\bibitem{WilFr}
\textbf{J Franks}, \textbf{R\,F Williams},
  \href{http://links.jstor.org/sici?sici=0002-9947(198709)303:1%3C97:BATJP%3E2%
.0.CO%3B2-N} {\emph{Braids and the {J}ones polynomial}}, Trans. Amer. Math.
  Soc. 303 (1987) 97--108 \xox{MR}{896009}

\bibitem{Gabai}
\textbf{D Gabai}, \href{http://dx.doi.org/10.1016/0040-9383(84)90001-6}
  {\emph{Foliations and genera of links}}, Topology 23 (1984) 381--394
  \xox{MR}{780731}

\bibitem{GLM}
\textbf{C\,M Gordon}, \textbf{R\,A Litherland}, \textbf{K Murasugi},
  \emph{Signatures of covering links}, Canad. J. Math. 33 (1981) 381--394
  \xox{MR}{617628}

\bibitem{Hirasawa}
\textbf{M Hirasawa}, \emph{The flat genus of links}, Kobe J. Math. 12 (1995)
  155--159 \xox{MR}{1391192}

\bibitem{Hirzebruch}
\textbf{F Hirzebruch}, \emph{Singularities and exotic spheres}, from:
  ``S\'eminaire Bourbaki, Vol.\ 10'', Soc. Math. France, Paris (1995)  Exp.\
  No.\ 314, 13--32 \xox{MR}{1610436}

\bibitem{Kawamura}
\textbf{T Kawamura}, \href{http://dx.doi.org/10.1007/s00014-002-8333-3}
  {\emph{Relations among the lowest degree of the {J}ones polynomial and
  geometric invariants for a closed positive braid}}, Comment. Math. Helv. 77
  (2002) 125--132 \xox{MR}{1898395}

\bibitem{Kreimer}
\textbf{D Kreimer}, \emph{Knots and {F}eynman diagrams}, Cambridge Lecture
  Notes in Physics 13, Cambridge University Press (2000) \xox{MR}{1778151}

\bibitem{KroMro}
\textbf{P\,B Kronheimer}, \textbf{T\,S Mrowka}, \emph{The genus of embedded
  surfaces in the projective plane}, Math. Res. Lett. 1 (1994) 797--808
  \xox{MR}{1306022}

\bibitem{Levine}
\textbf{J Levine}, \emph{Knot cobordism groups in codimension two}, Comment.
  Math. Helv. 44 (1969) 229--244 \xox{MR}{0246314}

\bibitem{MenThis}
\textbf{W\,W Menasco}, \textbf{M\,B Thistlethwaite}, \emph{The {T}ait flyping
  conjecture}, Bull. Amer. Math. Soc. $($N.S.$)$ 25 (1991) 403--412
  \xox{MR}{1098346}

\bibitem{Murasugi4}
\textbf{K Murasugi},
  \href{http://links.jstor.org/sici?sici=0003-486X(195903)2:69:2%3C258:GOALT%3%
E2.0.CO%3B2-B} {\emph{On the genus of the alternating knot. {I}, {II}}}, J.
  Math. Soc. Japan 10 (1958) 94--105, 235--248 \xox{MR}{0099664}

\bibitem{Murasugi3}
\textbf{K Murasugi},
  \href{http://links.jstor.org/sici?sici=0002-9327(196310)85:4%3C544:OACSOT%3E%
2.0.CO%3B2-S} {\emph{On a certain subgroup of the group of an alternating
  link}}, Amer. J. Math. 85 (1963) 544--550 \xox{MR}{0157375}

\bibitem{Murasugi2}
\textbf{K Murasugi},
  \href{http://links.jstor.org/sici?sici=0002-9947(196505)117%3C387:OACNIO%3E2%
.0.CO%3B2-P} {\emph{On a certain numerical invariant of link types}}, Trans.
  Amer. Math. Soc. 117 (1965) 387--422 \xox{MR}{0171275}

\bibitem{Nakamura}
\textbf{T Nakamura}, \href{http://dx.doi.org/10.1142/S0218216500000050}
  {\emph{Positive alternating links are positively alternating}}, J. Knot
  Theory Ramifications 9 (2000) 107--112 \xox{MR}{1749503}

\bibitem{Ozawa}
\textbf{M Ozawa}, \href{http://dx.doi.org/10.1007/s00014-002-8338-y}
  {\emph{Closed incompressible surfaces in the complements of positive knots}},
  Comment. Math. Helv. 77 (2002) 235--243 \xox{MR}{1915040}

\bibitem{OS}
\textbf{P Ozsv{\'a}th}, \textbf{Z Szab{\'o}},
  \href{http://dx.doi.org/10.2140/gt.2003.7.615} {\emph{Knot {F}loer homology
  and the four-ball genus}}, Geom. Topol. 7 (2003) 615--639 \xox{MR}{2026543}

\bibitem{Rasmussen}
\textbf{J Rasmussen}, \emph{Khovanov homology and the slice genus}
  \xox{arXiv}{math.GT/0402131}

\bibitem{Rolfsen}
\textbf{D Rolfsen}, \emph{Knots and links}, Mathematics Lecture Series 7,
  Publish or Perish, Berkeley, CA (1976) \xox{MR}{0515288}

\bibitem{Rudolph}
\textbf{L Rudolph}, \href{http://dx.doi.org/10.1016/0040-9383(82)90014-3}
  {\emph{Nontrivial positive braids have positive signature}}, Topology 21
  (1982) 325--327 \xox{MR}{649763}

\bibitem{Rudolph3}
\textbf{L Rudolph}, \emph{Quasipositivity as an obstruction to sliceness},
  Bull. Amer. Math. Soc. $($N.S.$)$ 29 (1993) 51--59 \xox{MR}{1193540}

\bibitem{Rudolph2}
\textbf{L Rudolph}, \emph{Positive links are strongly quasipositive}, from:
  ``Proceedings of the Kirbyfest (Berkeley, CA, 1998)'', (J Hass, M
  Scharlemann, editors), Geom. Topol. Monogr. 2 (1999)  555--562
  \xox{MR}{1734423}

\bibitem{gsig}
\textbf{A Stoimenow}, \emph{Bennequin's inequality and the positivity of the
  signature}, accepted by Trans. Amer. Math. Soc.

\bibitem{gener}
\textbf{A Stoimenow}, \emph{Diagram genus, generators and applications},
  preprint

\bibitem{gen2}
\textbf{A Stoimenow}, \emph{Knots of genus two} \xox{arXiv}{math.GT/0303012}

\bibitem{adeq}
\textbf{A Stoimenow}, \emph{On the crossing number of semiadequate links},
  preprint

\bibitem{2apos}
\textbf{A Stoimenow}, \href{http://dx.doi.org/10.1142/S0218216500000463}
  {\emph{The signature of 2-almost positive knots}}, J. Knot Theory
  Ramifications 9 (2000) 813--845 \xox{MR}{1775388}

\bibitem{gen1}
\textbf{A Stoimenow}, \href{http://dx.doi.org/10.1090/S0002-9939-01-05823-3}
  {\emph{Knots of genus one or on the number of alternating knots of given
  genus}}, Proc. Amer. Math. Soc. 129 (2001) 2141--2156 \xox{MR}{1825928}

\bibitem{posex_bcr}
\textbf{A Stoimenow}, \href{http://dx.doi.org/10.1090/S0002-9947-02-03022-2}
  {\emph{On the crossing number of positive knots and braids and braid index
  criteria of {J}ones and {M}orton-{W}illiams-{F}ranks}}, Trans. Amer. Math.
  Soc. 354 (2002) 3927--3954 \xox{MR}{1926860}

\bibitem{brlen}
\textbf{A Stoimenow}, \emph{The crossing number and maximal bridge length of a
  knot diagram}, Pacific J. Math. 210 (2003) 189--199 \xox{MR}{1989075}\ \ With
  an appendix by Mark Kidwell

\bibitem{pos}
\textbf{A Stoimenow}, \emph{Positive knots, closed braids and the {J}ones
  polynomial}, Ann. Sc. Norm. Super. Pisa Cl. Sci. $(5)$ 2 (2003) 237--285
  \xox{MR}{2004964}

\bibitem{apos}
\textbf{A Stoimenow}, \href{http://dx.doi.org/10.1112/S0010437X03000174}
  {\emph{Gau\ss\ diagram sums on almost positive knots}}, Compos. Math. 140
  (2004) 228--254 \xox{MR}{2004131}

\bibitem{STV}
\textbf{A Stoimenow}, \textbf{V Tchernov}, \textbf{A Vdovina},
  \href{http://dx.doi.org/10.1023/A:1021211008278} {\emph{The canonical genus
  of a classical and virtual knot}}, Geom. Dedicata 95 (2002) 215--225
  \xox{MR}{1950891}

\bibitem{SV}
\textbf{A Stoimenow}, \textbf{A Vdovina},
  \href{http://dx.doi.org/10.1007/s00208-005-0659-x} {\emph{Counting
  alternating knots by genus}}, Math. Ann. 333 (2005) 1--27 \xox{MR}{2169826}

\bibitem{Taniyama}
\textbf{K Taniyama}, \emph{A partial order of knots}, Tokyo J. Math. 12 (1989)
  205--229 \xox{MR}{1001742}

\bibitem{Traczyk}
\textbf{P Traczyk}, \href{http://dx.doi.org/10.1007/BF01258439}
  {\emph{Nontrivial negative links have positive signature}}, Manuscripta Math.
  61 (1988) 279--284 \xox{MR}{949818}

\bibitem{Vogel}
\textbf{P Vogel}, \emph{Representation of links by braids: a new algorithm},
  Comment. Math. Helv. 65 (1990) 104--113 \xox{MR}{1036132}

\bibitem{Williams}
\textbf{R\,F Williams}, \emph{Lorenz knots are prime}, Ergodic Theory Dynam.
  Systems 4 (1984) 147--163 \xox{MR}{758900}

\bibitem{Yamada}
\textbf{S Yamada}, \href{http://dx.doi.org/10.1007/BF01389082} {\emph{The
  minimal number of {S}eifert circles equals the braid index of a link}},
  Invent. Math. 89 (1987) 347--356 \xox{MR}{894383}

\bibitem{Yokota}
\textbf{Y Yokota}, \href{http://dx.doi.org/10.1016/0040-9383(92)90011-6}
  {\emph{Polynomial invariants of positive links}}, Topology 31 (1992) 805--811
  \xox{MR}{1191382}

\end{thebibliography}

\end{document}